\documentclass[12pt]{amsart}

\usepackage{color} % standard color package
\usepackage{graphicx} % standard graphics package
\usepackage{latexsym,amssymb,amsmath,amsfonts,bm} % extended symbol packages
\usepackage{mathrsfs} % mathscripts \mathscr
\usepackage{fixltx2e} % improves figure floating in Latex2e kernel
\usepackage{todonotes}
\usepackage{colonequals}

\allowdisplaybreaks[4]
\newtheorem{theorem}{Theorem}[section]
\newtheorem{lemma}{Lemma}[section]
\newtheorem{corollary}{Corollary}[section]
\newtheorem{propose}{Proposition}[section]
\newtheorem{definition}{Definition}[section]
\newtheorem{rem}{Remark}[section]
\newtheorem{example}{Example}[section]

\def\abs#1{\left\vert #1 \right\vert}
\def\allpolyell{\mbox{$\K^{\ell}\langle X \rangle$}}

\def\allpolyx0degn{\mbox{$P_n$}}

\def\allseries{\mbox{$\K\langle\langle X \rangle\rangle$}}
\def\allseriesR{\mbox{$\re\langle\langle X \rangle\rangle$}}
\def\allseriesC{\mbox{$\C\langle\langle X \rangle\rangle$}}
\def\allseriesdelta{\mbox{$\delta+\K\langle\langle X \rangle\rangle$}}
\def\allseriesdeltam{\mbox{$\delta+\K^m\langle\langle X \rangle\rangle$}}

\def\allseriesdeltaLC{\mbox{$\delta+\K_{LC}\langle\langle X \rangle\rangle$}}

\def\allseriesdeltamLC{\mbox{$\delta+\K^m_{LC}\langle\langle X \rangle\rangle$}}
\def\allseriesell{\mbox{$\K^{\ell} \langle\langle X \rangle\rangle$}}

\def\allseriesLC{\mbox{$\K_{LC}\langle\langle X \rangle\rangle$}}

\def\allseriesLCR{\mbox{$\re_{LC}\langle\langle X \rangle\rangle$}}
\def\allseriesLCC{\mbox{$\C_{LC}\langle\langle X \rangle\rangle$}}
\def\allseriesm{\mbox{$\K^m\langle\langle X \rangle\rangle$}}

\def\allseriesmLC{\mbox{$\K^{m}_{LC}\langle\langle X \rangle\rangle$}}

\def\allseriesellLC{\mbox{$\K^{\ell}_{LC}\langle\langle X \rangle\rangle$}}

\def\bfem#1{{\bf \em #1}} % boldface italics

\def\bull{\rule{0.08in}{0.08in}} % square filled bullet

\def\C{{\mathbb C}} % complex numbers (AMS symbol)
\newcommand{\comment}[1]{} % allows one to comment out a block of text

\def\Endallseries{{\rm End}(\allseries)}

\def\eqref#1{(\ref{#1})} % parentheses around referenced equation numbers
\def\Fliessdelta{\mathscr{F}_{\delta}}

\def\id{{\rm id}}

\def\mbf#1{\hbox{\mathversion{bold}$#1$}} % math boldface
\def\modcomp{\:\tilde{\circ}\,} % modified composition product

\def\nat{{\mathbb N}} % natural numbers (AMS symbol)
\def\norm#1{\left\Vert#1\right\Vert}

\def\openbull{\framebox[0.08in][c]{$\;$}} % square unfilled bullet

\def\re{{\mathbb R}} % real numbers (AMS symbol)

\def\shuffle{{\scriptscriptstyle \;\sqcup \hspace*{-0.05cm}\sqcup\;}}

\def\supp{{\rm supp}}
\def\begals{\[\begin{aligned}}
\def\endals{\end{aligned}\]}
\def\begce{\begin{center}}
\def\endce{\end{center}}
\def\begar{\begin{array}}
\def\endar{\end{array}}
\def\begeq{\begin{equation}}
\def\endeq{\end{equation}}
\def\begdi{\begin{displaymath}}
\def\enddi{\end{displaymath}}
\def\begdis{\begin{eqnarray*}}
\def\enddis{\end{eqnarray*}}
\def\begeqa{\begin{eqnarray}}
\def\endeqa{\end{eqnarray}}
\def\begdes{\begin{description}}
\def\enddes{\end{description}}
\def\begit{\begin{itemize}}
\def\endit{\end{itemize}}
\def\begen{\begin{enumerate}}
\def\enden{\end{enumerate}}
\def\beglar{\left[\begin{array}}
\def\endrar{\end{array}\right]}
\def\begle{\begin{lemma}}
\def\endle{\end{lemma}}
\def\begde{\begin{definition}}
\def\endde{\end{definition}}
\def\begth{\begin{theorem}}
\def\endth{\end{theorem}}
\def\begco{\begin{corollary}}
\def\endco{\end{corollary}}
\def\begprop{\begin{prop}}
\def\endprop{\end{prop}}
\def\begex{\begin{example}}
\def\endex{\hfill\openbull \end{example} \vspace*{0.15in}}
\def\begexer{\begin{myexercise}}
\def\endexer{\end{myexercise}}

\def\begres{\noindent{\bf Remarks}:\begin{enumerate}}
\def\endres{\end{enumerate} \par}
\def\begpr{\noindent{\em Proof:}$\;\;$}
\def\endpr{\hfill\bull \vspace*{0.15in}}
\def\begtab{\begin{tabular}}
\def\endtab{\end{tabular}}
\def\rref#1{(\ref{#1})}
 \newcommand{\ellInfty}[3]{\ell_{\infty,#1}(#2,#3)}
 \newcommand{\ellInftyNorm}[2]{\norm{#1}_{\ell_{\infty,#2}}}
 
 \newcommand{\R}{\mathbb{R}}
 \newcommand{\Frechet}{{Fr\'echet} }
\newcommand{\LB}[1][\cdot \hspace{1pt} , \cdot]{[\hspace{1pt} #1 \hspace{1pt} ]}

\DeclareMathOperator{\evol}{evol}
\DeclareMathOperator{\Evol}{Evol}
\newcommand{\K}{\ensuremath{\mathbb{K}}}
\DeclareMathOperator{\Lf}{\mathbf{L}}
\def\shuff#1#2{\mathbin{
      \hbox{\vbox{\hbox{\vrule \hskip#2 \vrule height#1 width 0pt}\hrule}\vbox{\hbox{\vrule \hskip#2 \vrule height#1 width 0pt\vrule }\hrule}}}}
\def\shuffl{{\mathchoice{\shuff{5pt}{3.5pt}}{\shuff{5pt}{3.5pt}}{\shuff{3pt}{2.6pt}}{\shuff{3pt}{2.6pt}}}}
\def\shuffle{{\, \shuffl \,}}

\newcommand{\coloneq}{\colonequals}

\begin{document}

\title{Continuity of Formal Power Series Products in Nonlinear Control Theory}

\author{W. Steven Gray}
\address{Old Dominion University, Norfolk, Virginia 23529, USA}
\email{sgray@odu.edu}

\author{Mathias Palmstr{\o}m}
\address{Department of Mathematics, Universitet i Bergen, All\'{e}gate 41, 5020 Bergen, Norway}
\email{Mathias.Palmstrom@student.uib.no}

\author{Alexander Schmeding}
\address{FLU, Nord university, H\o{}gskoleveien 27, 7601 Levanger, Norway}
\email{alexander.schmeding@nord.no}

\date{}

\begin{abstract} % Abstract of not more than 250 words.
Formal power series products appear in nonlinear control theory when systems modeled by Chen-Fliess series are interconnected to form new systems. In fields like adaptive control and learning systems, the coefficients of these formal power series are estimated sequentially with real-time data. The main goal is to prove the
continuity and analyticity of such products with respect to several natural (locally convex) topologies on spaces of locally convergent formal power series in order to
establish foundational properties behind these technologies. In addition, it is shown that a transformation group central to describing the output feedback connection is
in fact an analytic Lie group in this setting with certain regularity properties.
\end{abstract}

\maketitle
\textbf{MSC2020}: 93C10 (primary), % Nonlinear systems in control theory
46A04, % Locally convex Fr´echet spaces and (DF)-spaces
46A13, % Spaces defined by inductive or projective limits (LB, LF, etc.)
47N70, % Applications of operator theory in systems, signals, circuits, and control theory
22E65, % Infinite-dimensional Lie groups
46B45, % Banach sequence spaces
16T30  % Hopf algebras and combinatorics
\smallskip

\textbf{Keywords}: nonlinear control systems, Chen--Fliess series, system interconnection, Silva space, real analytic, locally convex Lie group, regularity of Lie groups

\noindent
\tableofcontents

\thispagestyle{empty}

\section{Introduction}

The interconnection of simple input-output systems to form more complex and useful systems is commonplace in science and engineering.
When each component is a nonlinear dynamical system, a weighted infinite sum of iterated integrals known as a
{\em Chen-Fliess series}, $F_c$, provides a convenient way to represent its local behavior \cite{Fliess_81,Fliess_83,Fliess-et.al._83,Lamnabhi-Lagarrigue_96,Wang_90}.
When this series converges on some set of admissible inputs, $F_c$ defines a so called {\em Fliess operator}.
It is uniquely specified by a formal power series $c$ in $\allseries$, known as its {\em generating series},
where $X$ is a finite set of indeterminants, and $\K$ is a suitable field.
An interconnection of two Chen-Fliess series $F_c$ and $F_d$, represented by $F_c\Box F_d$, induces a corresponding algebra $(\allseries,\Box^\prime)$
so that $F_c\Box F_d=F_{c\Box^\prime d}$ \cite{Ferfera_79,Ferfera_80,Fliess_81,Gray-etal_SCL14,Gray-Li_05}.
Algebras defined in this manner provide computational frameworks for explicitly computing the generating series of interconnected
systems for the purposes of analysis and design, especially in the field of nonlinear control theory.
Historically, the coefficients of $c$ have been determined by direct calculations using state space models
derived from physical laws and other first principles \cite{Isidori_95,Nijmeijer-vanderSchaft_90}. But with the growth of adaptive control and new types of learning
based technologies, there is increasing interest in estimating these coefficients using real-time data and numerical methods from the field of system identification
\cite{Gray-etal_Automatica_20,Padoan-Astolfi_16}. Assuming that a given sequence of estimates
asymptotically approaches its true value as more data is collected, a difficult problem in
its own right, there is a fundamental question regarding continuity. Consider a sequence of generating series $c_i$, $i>1$
known to produce a sequence of corresponding Fliess operators $F_{c_i}$, $i\geq 1$.
If $c_i\rightarrow c$ in some manner, is it also true that $F_{c_i}$, $i\geq 1$ converges to a well defined Fliess operator $F_c$, i.e., does
the limit point $c$ ensure a convergent Chen-Fliess series?
The answer, of course, depends directly on the ambient sets and the assumed topologies.
For example, in \cite{Dahmen-etal_MTNS21,Winter_Arboleda_19} the claim is shown to be false on the subset of {\em locally convergent} series in $\allseries$, (i.e.,
a set of generating series under which their corresponding Fliess operators are known to converge at least locally)
endowed with the ultrametric topology and where the operator
space has an $L_p$ type topology. As the ultrametric topology mirrors the algebra but provides almost no information on the analytic behavior of the series, this
outcome is not surprising.
On the other hand, in \cite{Dahmen-etal_MTNS21} the claim is shown to be true
when the ultrametric topology is replaced with a certain Banach topology on a subspace. Ultimately, the question boils down to
identifying topological vectors spaces contained in $\allseries$ which ensure that every limit point is a generating series with a well defined Fliess operator in some
sense.

The main goal of this paper is to address a natural follow-up question: Suppose $c_i,d_i\in\allseries$, $i\geq 1$ are two sequences of generating series
converging to $c$ and $d$, respectively, such that $F_{c_i}$ and $F_{d_i}$, $i\geq 1$ are well defined Fliess operators as are its limit points.
Assuming that $c\Box^\prime d$ has the same convergence properties as $c$ and $d$ (this theory is well understood, see \cite{Gray-Wang_SCL02,Thitsa-Gray_SIAM12,Winter_Arboleda_19}),
under what conditions does $c_i\Box^\prime d_i\mapsto c\Box^\prime d$? Is it even possible to identify infinite-dimensional spaces with respect to which the products are smooth or even analytic?

Three formal power series products will be considered: the shuffle product, which models a type of parallel connection \cite{Fliess_81};
a composition product modeling series connections \cite{Ferfera_79,Ferfera_80,Gray-Li_05}; and a group product for a transformation group
known to model dynamic output feedback, a central object of study in control theory \cite{Gray-etal_SCL14}. In addition, the continuity of the shuffle inverse will be addressed.
(A preliminary version of this analysis was presented in \cite{palm20}.)
The shuffle group appears in the context of feedback linearization \cite{Gray-etal_AUTO14,Gray-Ebrahimi-Fard_SIAM17}.
In each case continuity will be considered in both the Fr\'{e}chet and Silva topologies. In addition, analyticity of these
product will be characterized. It should be noted that the Fr\'{e}chet topology was used in \cite{Winter_Arboleda_19} to show that the shuffle and composition products preserve
a type of {\em global convergence}. Continuity issues in this setting are beyond the scope of the present paper. However, the \Frechet topology is employed as a natural (locally convex) topology on the space of all power series. Convergence in this topology does not preserve growth bounds. Thus, it is necessary to endow the space of locally convergent series with the finer Silva topology.
Next it will be shown that the output feedback transformation group is a {\em locally convex} Lie group (see \cite{neeb2006} for a survey on (infinite-dimensional) Lie theory).
This result builds on the development of a pre-Lie algebra presented in \cite{DuffautEspinosa2016609,Foissy_15}.
Lie groups have a long history in feedback control theory originating with the work of Brockett in \cite{Brockett_76}.
More recent applications in this context have appeared in \cite{Gray-Ebrahimi-Fard_SIAM17,Gray-Ebrahimi-Fard_SCL21}, albeit only in the
formal case where an explicit differential structure is not specified. The present work will provide a means to fill this gap.
Finally, the regularity of these Lie groups is investigated. Roughly speaking, regularity of a Lie group asks for the existence and smooth parameter dependence of certain ordinary differential equations on the Lie group. Note that since the Lie groups at hand are not modeled on Banach spaces, the usual theory for existence and uniqueness of ordinary differential equations does not apply. However, it is shown that the \Frechet Lie groups are regular. While some progress on the regularity problem is made for the Silva Lie groups, their regularity largely remains an open problem that the authors plan to pursue in future work.

\textbf{Acknowledgements} A.S.~wishes to thank the University of Bergen, Norway, where he was employed while most of the present work was carried out.
\section{Preliminaries}

Throughout this paper let $\K \in \{\R,\C\}$, namely either the field of real numbers $\R$ or the field of complex numbers $\C$.
It will be essential to admit complex coefficients in order to discuss analyticity of mappings on infinite-dimensional spaces. Note that the continuity results are unaffected by this choice. Refer to Appendix \ref{app:smooth} for more information regarding calculus on infinite-dimensional spaces.

\subsection{Chen-Fliess series}

An \emph{alphabet} $X=\{ x_0,x_1,$ $\ldots,x_m\}$ is any nonempty and finite set
of noncommuting symbols referred to as \emph{letters}. A \emph{word} $\eta=x_{i_1}\cdots x_{i_k}$ is a finite sequence of letters from $X$.
The number of letters in a word $\eta$, written as $\abs{\eta}$, is called its \emph{length}.
The empty word, $\emptyset$, is taken to have length zero.
The collection of all words having length $k$ is denoted by
$X^k$. Define the set of all words $X^\ast=\bigcup_{k\geq 0} X^k$, which constitutes a monoid under the concatenation product.
Any mapping $c\colon X^\ast\rightarrow
\K^\ell$ is called a \emph{formal power series}.
Often $c$ is
written as the formal sum $c=\sum_{\eta\in X^\ast}(c,\eta)\eta$,
where the \emph{coefficient} $(c,\eta)$ is the image of
$\eta\in X^\ast$ under $c$.
The \emph{support} of $c$, $\supp(c)$, is the set of all words having nonzero coefficients.
A series $c$ is said to be {\em proper} when $\emptyset\not\in \supp(c)$.
The set of all noncommutative formal power series over the alphabet $X$ is
denoted by $\allseriesell$. The subset of series with finite support, i.e., polynomials,
is represented by $\allpolyell$.
Each set is an associative $\K$-algebra under the catenation product and an associative and commutative $\K$-algebra under the \emph{shuffle product}, that is, the bilinear product uniquely specified by the shuffle product of two words
\begdi
(x_i\eta)\shuffle(x_j\xi)=x_i(\eta\shuffle(x_j\xi))+x_j((x_i\eta)\shuffle \xi),
\enddi
where $x_i,x_j\in X$, $\eta,\xi\in X^\ast$ and with $\eta\shuffle\emptyset=\emptyset\shuffle\eta=\eta$ \cite{Fliess_81}.
On $\allseriesell$ the definition is extended componentwise.

$\allseriesell$ can be viewed as a locally convex space whose topology is briefly described next. First note that identifying a formal power series with the sequence of its coefficients
defines an isomorphism of vector spaces $\allseriesell \cong \prod_{\eta \in X^\ast}\K^\ell$. The space on the right hand side is a countable product of Banach spaces, hence a complete metrisable locally convex vector space (i.e., a \Frechet space). Thus, $\allseriesell$ inherits a canonical \Frechet space structure. By construction the evaluations $a_\eta \colon \allseriesell \rightarrow \K^\ell, c \mapsto (c,\eta)$ are continuous. Therefore, convergence in this topology is equivalent to separate convergence of all coefficients of a series towards the corresponding coefficients of the limit series. Moreover, the \Frechet topology is initial with respect to the point evaluations, i.e., a map $f$ to $\allseriesell$ is continuous if and only if $a_{\eta} \circ f$ is continuous for every word $\eta \in X^\ast$.

Given any $c\in\allseriesell$ one can associate a causal
$m$-input, $\ell$-output operator, $F_c$, in the following manner.
Let $\mathfrak{p}\ge 1$ and $t_0 < t_1$ be given. For a Lebesgue measurable
function $u: [t_0,t_1] \rightarrow\K^m$, define
$\norm{u}_{\mathfrak{p}}=\max\{\norm{u_i}_{\mathfrak{p}}: \ 1\le
i\le m\}$, where $\norm{u_i}_{\mathfrak{p}}$ is the usual
$L_{\mathfrak{p}}$-norm for a measurable real-valued function,
$u_i$, defined on $[t_0,t_1]$.  Let $L^m_{\mathfrak{p}}[t_0,t_1]$
denote the set of all measurable functions defined on $[t_0,t_1]$
having a finite $\norm{\cdot}_{\mathfrak{p}}$ norm and
$B_{\mathfrak{p}}^m(R_u)[t_0,t_1]:=\{u\in
L_{\mathfrak{p}}^m[t_0,t_1]:\norm{u}_{\mathfrak{p}}\leq R_u\}$.
Assume $C[t_0,t_1]$
is the subset of continuous functions in $L_{1}^m[t_0,t_1]$. Define
inductively for each $\eta\in X^{\ast}$ the map $E_\eta:
L_1^m[t_0, t_1]\rightarrow C[t_0, t_1]$ by setting
$E_\emptyset[u]=1$ and letting
\[E_{x_i\bar{\eta}}[u](t,t_0) =
\int_{t_0}^tu_{i}(\tau)E_{\bar{\eta}}[u](\tau,t_0)\,d\tau, \] where
$x_i\in X$, $\bar{\eta}\in X^{\ast}$, and $u_0=1$. The
\emph{Chen-Fliess series} corresponding to $c$ is
\begeq
y(t)=F_c[u](t) =
\sum_{\eta\in X^{\ast}} (c,\eta)\,E_\eta[u](t,t_0) \label{eq:Fliess-operator-defined}
\endeq
\cite{Fliess_81,Fliess_83}.
It can be shown that if there exists real numbers $K,M\geq 0$ such that
\begeq \label{eq:LC-condition}
\abs{(c,\eta)}\leq KM^{\abs{\eta}}\abs{\eta}!,\;\;\forall \eta\in X^\ast
\endeq
($\abs{z}:=\max_i \abs{z_i}$ when $z\in\K^\ell$)
then the series defining $F_c$ converges absolutely and uniformly for sufficient
small $R,T >0$ and
constitutes a well defined mapping from
$B_1^m(R)[t_0,$ $t_0+T]$ into $B_{\infty}^{\ell}(S)[t_0, \, t_0+T]$ for some $S>0$.
Any such mapping is called a \emph{locally convergent Fliess operator}.
Here $\allseriesellLC$ will denote the set of all such
{\em locally convergent} generating series, i.e., those series
satisfying growth condition \rref{eq:LC-condition}.
Given any smooth state space realization of $y=F_c[u]$,
\begdi
\dot{z}=g_0(z)+\sum_{i=1}^m g_i(z)u_i,\;\;z(0)=z_0,\;\;y=h(z),
\enddi
it is known that the generating series $c$ is determined by
\begeq \label{eq:c-from-Lgh}
(c_j,\eta)=L_{g_{i_1}}\cdots L_{g_{i_k}}h_j(z_0),\;\;\eta=x_{i_k}\cdots x_{i_1}\in X^\ast,\;\;j=1,2,\ldots,\ell
\endeq
where $L_{g_i}h_j$ is the Lie derivative of $h_j$ with respect to $g_i$.

\subsection{Formal power series products induced by system interconnection}
\label{subsec:system-interconnections}

Given Fliess operators $F_c$ and $F_d$, where $c,d\in\allseriesellLC$,
the parallel and product connections satisfy $F_c+F_d=F_{c+d}$ and $F_cF_d=F_{c\shuffle d}$,
respectively \cite{Fliess_81}.
When Fliess operators $F_c$ and $F_d$ with $c\in\allseriesellLC$ and
$d\in\allseriesmLC$ are interconnected in a cascade fashion, the composite
system $F_c\circ F_d$ has the
Fliess operator representation $F_{c\circ d}$, where
the {\em composition product} of $c$ and $d$
is given by
\begeq \label{eq:c-circ-d}
c\circ d=\sum_{\eta\in X^\ast} (c,\eta)\,\psi_d(\eta)(\mbf{1})
\endeq%
\cite{Ferfera_79,Ferfera_80}. Here $\mbf{1}$ denotes the monomial $1\emptyset$, and
$\psi_d$ is the continuous (in the ultrametric sense) algebra homomorphism
from $\allseries$ to the set of vector space endomorphisms on $\allseries$, $\Endallseries$, uniquely specified by
$\psi_d(x_i\eta)=\psi_d(x_i)\circ \psi_d(\eta)$ with
$ %\label{eq:psi-d-on-words}
\psi_d(x_i)(e)=x_0(d[i]\shuffle e),
$
$i=0,1,\ldots,m$
for any $e\in\allseries$,
and where $d[i]$ is the $i$-th component series of $d$
($d[0]:=\mbf{1}$). By definition,
$\psi_d(\emptyset)$ is the identity map on $\allseries$.

When two Fliess operators $F_c$ and $F_d$ are interconnected
to form a feedback system with $F_c$ in the forward path and $F_d$ in the
feedback path, the generating series of the closed-loop system is
denoted by the {\em feedback product} $c@d$. It can
be computed explicitly using the Hopf algebra of coordinate functions
associated with the underlying {\em output feedback group} \cite{Gray-etal_SCL14}.
Specifically, in the single-input, single-output case where $X=\{x_0,x_1\}$ and $\ell=1$, define the set of
{\em unital} Fliess operators
$
\Fliessdelta=\{I+F_c\;:\;c\in\allseriesLC\},
$
where $I$ denotes the identity map.
It is convenient to introduce the symbol
$\delta$ as the (fictitious) generating series for the identity map. That is,
$F_\delta:=I$ such that $I+F_c:=F_{\delta+c}=F_{c_\delta}$ with
$c_\delta:=\delta+c$.
The set of all such generating series for
$\Fliessdelta$ will be denoted by $\allseriesdeltaLC$.
The central idea is that $(\Fliessdelta,\circ,I)$ forms a group of operators
under the composition
\begdi
F_{c_\delta}\circ F_{d_\delta}=(I+F_c)\circ(I+F_d)
= F_{c_\delta\circ d_\delta},
\enddi
where $c_\delta\circ d_\delta:=\delta+c\circledcirc d$, $c\circledcirc d:=d+c\modcomp d_\delta$, and
$\modcomp$ denotes the {\em mixed} composition product.
That is, the product
\begeq \label{eq:mixed-composition-product}
c\modcomp d_\delta =\sum_{\eta\in X^\ast} (c,\eta)\,\phi_d(\eta)(\mbf{1}),
\endeq
where $\phi_d$ is analogous to $\psi_d$ in \rref{eq:c-circ-d}
except here $\phi_d(x_i)(e)=x_ie+x_0(d[i]\shuffle e)$ with $d[0]:=0$ \cite{Gray-Li_05}.
The set of unital generating series $\allseriesdelta$ (not necessarily locally convergent) forms
a group $(\allseriesdelta,\circ,\delta)$. The restriction to the set of locally convergent series defines the
subgroup $\allseriesdeltaLC$. The
mixed composition product can be viewed as a right action of $\allseriesdelta$ acting freely on $\allseries$ \cite{Gray-etal_CDC13}.
The corresponding Hopf algebra $H$ is the free algebra generated by
the coordinate maps
\begdi \label{eq:character-maps}
a_\eta:\allseriesdelta\rightarrow \K:c_\delta\mapsto (c,\eta),\quad \eta\in X^\ast
\enddi
under the commutative product
\begdi \label{eq:mu-product}
\mu:a_\eta\otimes a_{\xi}\mapsto a_{\eta}a_{\xi},
\enddi
where the unit $\mbf{1}_\delta$ is defined to map every $c_\delta$ to one. Let $V$ be the $\K$-vector space
of coordinate functions.
If the {\em degree} of $a_{\eta}$ is defined as
$\deg(a_{\eta})=2\abs{\eta}_{x_0}+\abs{\eta}_{x_1}+1$, then both $V$ and the
algebra $H$ are graded and connected with $V=\bigoplus_{n\geq 0} V_n$ and
$H=\bigoplus_{n\geq 0} H_n$, where $V_n$ and $H_n$ are sets containing all the
degree $n$ elements, and $V_0=H_0=\K\mbf{1}_\delta$.
The coproduct $\Delta$ is defined so that
\begdi
\Delta a_\eta(c_\delta,d_\delta)=a_\eta(c_\delta\circ d_\delta)=(c_\delta\circ d_\delta,\eta).
\enddi
Of primary importance is the following lemma which describes how
the group inverse $c_\delta^{\circ -1}:=\delta+c^{\circ -1}$ is computed.

\begin{lemma} \cite{Gray-etal_SCL14} \label{le:antipode-is-group-inverse}
The Hopf algebra $(H,\mu,\Delta)$ has an antipode $S$ satisfying
$a_\eta(c_\delta^{\circ -1})=(Sa_\eta)(c_\delta)$ for all $\eta\in X^\ast$ and $c_\delta\in\allseriesdelta$.
\end{lemma}

With this concept, the generating series for the feedback connection, $c@d$, can be computed
explicitly as described in the next theorem. It states that feedback in the present
context can be viewed in terms of the group $(\allseriesdelta,\circ,\delta)$ acting on $\allseries$
in a specific manner.

\begin{theorem} \cite{Gray-etal_SCL14}
\label{th:feedback-product-formula}
For any $c,d\in\allseries$ it follows that
\begdi
c@d=c\modcomp(-d\circ c)_\delta^{\circ -1}. \label{eq:catd-formula}
\enddi
\end{theorem}

In addition to the elementary system interconnections described above,
there is the quotient connection that is useful in the
context of system inversion \cite{Gray-etal_AUTO14}.
This is a type of parallel connection where the quotient of the subsystems' outputs
is computed. In terms of generating series, the quotient is realized using
the shuffle inverse as described next. Division by zero is avoided by requiring
the divisor series to be non proper.

\begin{theorem} \cite{Gray-etal_AUTO14} \label{th:shuffle-group}
The set of non proper series in $\allseries$ is a group under the shuffle
product. In particular, the shuffle inverse of any such series $c$ is
\begdi
c^{\shuffle -1}=((c,\emptyset)(1-c^\prime))^{\shuffle -1}=(c,\emptyset)^{-1}(c^{\prime})^{\shuffle\ast},
\enddi
where $c^\prime:=\mbf{1}-c/(c,\emptyset)$ is proper, and $(c^\prime)^{\shuffle\ast}:=\sum_{k\geq 0} (c^\prime)^{\shuffle k}$.
\end{theorem}

\begin{theorem} \cite{Gray-etal_AUTO14} \label{th:quotient-connection}
For $c,d\in\allseriesLC$, the quotient connection $F_c/F_d$ has a Fliess operator
representation if and only if $d$ is non proper. In particular, $F_c/F_d=F_{c/d}$,
where $c/d:=c\shuffle d^{\shuffle -1}$.
\end{theorem}

\section{Continuity of formal power series products}

In this section, the continuity of the various products modeling system interconnections described in the previous section is proved.
The main goal is to establish continuity on spaces of locally convergent series. In \cite{Dahmen-etal_MTNS21} the authors described the space of locally
convergent Chen-Fliess series as a locally convex space carrying a Silva space topology. That construction is summarized first, and then the continuity
results are presented.

Fix $M>0$ and define
\[
       \ellInftyNorm{c}{M} := \sup\left\{   \frac{\abs{(c,\eta)}}{M^{\abs{\eta}}\abs{\eta}!} : \eta\in X^\ast    \right\} \in \left[0,\infty\right]
\]
for each $c\in\allseriesell$.
The set of all $c$ with $\ellInftyNorm{c}{M}<\infty$ is denoted by $\ellInfty{M}{X^*}{\K^\ell}$.
It is straightforward to check that $\ellInfty{M}{X^*}{\K^\ell}$
is a vector subspace of $\allseriesell$. The function $\ellInftyNorm{\cdot}{M}$ is a norm on $\ellInfty{M}{X^*}{\K^\ell}$.
This space is a Banach space as it is isometrically isomorphic to the Banach space of all bounded functions
$\ell_\infty(X^*,\K^\ell):= \{c \colon X^* \rightarrow \K^\ell : \sup_{\eta }\abs{(c,\eta)}<\infty\}$.
The Banach space of generating series bounded with respect to the constant $M$ obviously does not capture all locally convergent series.
Indeed for larger $M$ one obtains series which converge only on a smaller disc. To capture all locally convergent series in one space, it is necessary to
pass to the limit of these Banach spaces as described next.

\begin{definition} [Locally convergent series as a Silva space]
Consider the union \[
  \ellInfty{\rightarrow}{X^*}{\K^\ell} := \bigcup_{M>0} \ellInfty{M}{X^*}{\K^\ell}.
\]
Topologise this space as the locally convex inductive limit of the system
$\left(\ellInfty{M}{X^*}{\K^\ell} \right)_{M>0}$.
\end{definition}

One can show that the inclusion mappings in this sequence are compact operators, hence the resulting space is a \emph{Silva space} \cite{BS16,DS18}.
Since the sequence $M_k = k$, $k\in \mathbb N$ is cofinal, one can always find an $M\in\mathbb N$ for which $\ellInftyNorm{c}{M}<\infty$.
Thus, one could equivalently work only with $M\in\mathbb N$.
Though the Silva space topology is more complicated than the Banach spaces from which it was built, some of its properties make it very amenable for the applications
considered here.
The most important properties are summarized in the next lemma. Refer to \cite{Yos57} for proofs and more information about Silva spaces.

\begin{lemma} [Properties of Silva spaces] \label{lem:Silvaprop}
\rule{0in}{0in}
\begin{enumerate}
 \item A sequence converges in $ \ellInfty{\rightarrow}{X^*}{\K^\ell}$ if and only if there exists $M>0$ such that the sequence is contained and converges in the Banach space $\ellInfty{M}{X^*}{\K^\ell}$.
 \item Silva spaces are sequential, thus a map defined on a Silva space is continuous if and only if it is sequentially continuous. Moreover, Silva spaces are separable and finite products of Silva spaces are again Silva spaces.
 \item A mapping $f \colon \ellInfty{\rightarrow}{X^*}{\K^\ell} \rightarrow E $ into a locally convex space is continuous (differentiable) if and only if for every $M>0$ the induced mapping $$f_M := f|_{\ellInfty{M}{X^*}{\K^\ell}} \colon\ellInfty{M}{X^*}{\K^\ell} \rightarrow E$$ is continuous (differentiable).
\end{enumerate}
\end{lemma}

Perhaps the most striking property of the Silva topology is that one can address continuity and differentiability questions in the Banach spaces from which the Silva space is built.
This will be demonstrated in the next section addressing the continuity of formal power series products.

\subsection{Continuity of shuffle product and shuffle inverse}

The following lemma is a prerequisite for proving continuity of the shuffle product.

\begin{lemma} \label{le:shuffle-Ke}
Fix $M>0$. If $c,d\in\ell_{\infty,M}(X^\ast,\K^\ell)$, then $c\shuffle d\in \ell_{\infty,M_\epsilon}(X^\ast,\K^\ell)$
for any $M_\epsilon=M(1+\epsilon)$, $\epsilon>0$
and
\begdis
\|c\shuffle d\|_{\ell_\infty,M_\epsilon}\leq K_\epsilon\,\|c\|_{\ell_\infty,M}\|d\|_{\ell_\infty,M},
\enddis
where
$K_\epsilon=\sup_{\eta\in X^\ast}({\abs{\eta}+1})/(1+\epsilon)^{\abs{\eta}}\leq \hat{K}_\epsilon:=e^{-1}(1+\epsilon)/(\log(1+\epsilon))$.
\end{lemma}

\begin{figure}[h]
\begin{center}
\includegraphics[width=7cm]{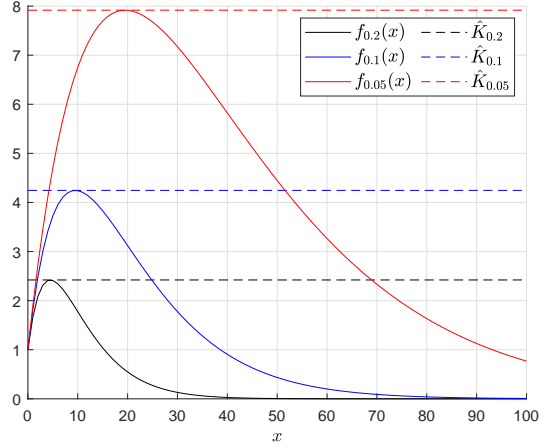}
\end{center}
\caption{Sample plots of $f_\epsilon(x)$ and $\hat{K}_\epsilon$ in Lemma~\ref{le:shuffle-Ke}}.
\label{fig:ke-plot}
\end{figure}

\begpr
For any $\eta\in X^\ast$
\begin{align*}
\abs{(c\shuffle d,\eta)}&=\abs{\sum_{k=0}^{|\eta|}
\sum_{\nu\in X^k \atop \xi\in X^{|\eta|-k}}
(c,\nu) (d,\xi)(\nu\shuffle\xi,\eta)} \\
&\leq \sum_{k=0}^{|\eta|}
\sum_{\nu\in X^k \atop \xi\in X^{|\eta|-k}}
 \|c\|_{\ell_\infty,M}M^k k!\;  \|d\|_{\ell_\infty,M}M^{\abs{\eta}-k} (\abs{\eta}-k)!\; (\nu\shuffle\xi,\eta) \\
&=  \|c\|_{\ell_\infty,M} \|d\|_{\ell_\infty,M}M^{\abs{\eta}} \sum_{k=0}^{\abs{\eta}} k!\:(\abs{\eta}-k)!\:{\abs{\eta}
\choose k} \\
&=\|c\|_{\ell_\infty,M} \|d\|_{\ell_\infty,M}M^{\abs{\eta}} \sum_{k=0}^{\eta} \abs{\eta} ! \\
&=\|c\|_{\ell_\infty,M} \|d\|_{\ell_\infty,M}M^{\abs{\eta}} (\abs{\eta}+1)!.
\end{align*}
Note that this bound is achievable when $c=\sum_{\eta\in X^\ast} K_c M^{|\eta|}|\eta|!\,\eta$ and
$d=\sum_{\eta\in X^\ast} K_d M^{|\eta|}|\eta|!\,\eta$ for any $K_c,K_d\geq 0$. Now define $M_\epsilon=M(1+\epsilon)$ with $\epsilon>0$
and rewrite the final inequality above as
\begin{align*}
\frac{\abs{(c\shuffle d,\eta)}}{M_\epsilon^{\abs{\eta}}\abs{\eta}!}&\leq \|c\|_{\ell_\infty,M} \|d\|_{\ell_\infty,M}\frac{{\abs{\eta}+1}}{(1+\epsilon)^{\abs{\eta}}},\;\;\forall \eta\in X^\ast.
\end{align*}
Taking the supremum over $X^\ast$ gives
\begin{align*}
\|c\shuffle d\|_{\ell_\infty,M_\epsilon}&\leq K_{\epsilon}\, \|c\|_{\ell_\infty,M} \|d\|_{\ell_\infty,M},\;\;\forall \eta\in X^\ast,
\end{align*}
where $K_\epsilon=\sup_{\eta\in X^\ast}({\abs{\eta}+1})/(1+\epsilon)^{\abs{\eta}}$. The upper bound for $K_\epsilon$ is found by showing that
$f_\epsilon(x)=({x+1})/(1+\epsilon)^{x}$ has a single maximum at $x_\epsilon^\ast=(1/\log(1+\epsilon))-1>0$ when $0<\epsilon\leq e-1$, and $\hat{K}_\epsilon=f_\epsilon(x_\epsilon^\ast)=e^{-1}(1+\epsilon)/(\log(1+\epsilon))$.
In this case, the upper bound is tight (see Figure~\ref{fig:ke-plot}). For $\epsilon > e-1$, $K_\epsilon=1$ and $\hat{K}_\epsilon >1$, and thus this upper bound
is conservative.
\endpr

\begin{theorem} \label{thm:shuffle:cont}
The shuffle product is continuous on $\allseries$ and $\allseriesLC$ with respect to the \Frechet and the Silva topology, respectively.
\end{theorem}

\begpr
Consider the shuffle product on $\allseries$. Since the topology of $\allseries$ is initial with respect to the coordinate functions $a_\eta \colon \allseries \rightarrow \K, c \mapsto (c,\eta)$, it suffices to prove that $a_\eta \circ \shuffle$ is continuous for each $\eta \in X^\ast$. However, as was seen in the proof of Lemma \ref{le:shuffle-Ke} for $\eta \in X^\ast$, it follows that
\begin{align*}
 a_\eta \circ \shuffle (c,d) &= (c\shuffle d,\eta) = \sum_{k=0}^{|\eta|} \sum_{\nu\in X^k \atop \xi\in X^{|\eta|-k}} (c,\nu) (d,\xi) (\nu \shuffle \xi, \eta)\\
 &=  \sum_{k=0}^{|\eta|} \sum_{\nu\in X^k \atop \xi\in X^{|\eta|-k}} a_\nu(c)a_\xi (d) (\nu \shuffle \xi, \eta). %\label{shuffle:poly}
\end{align*}
This  shows that $a_\eta \circ \shuffle (c,d)$ is a polynomial in the variables $a_\nu(c),a_\xi(d)$. Since the coordinate functions are continuous in the series $c,d$, it is clear that the shuffle product is continuous. Thus, the shuffle product is continuous on $\allseries$. For the corresponding result on $\allseriesLC$, apply Lemma~\ref{le:shuffle-Ke} for any $\epsilon>0$:
\begin{align*}
\|(c\shuffle d)-(c_j\shuffle d_j)\|_{\ell_\infty,M_\epsilon}&=\|(c-c_j)\shuffle d+c_j\shuffle (d-d_j)\|_{\ell_\infty,M_\epsilon}\\
&\leq \|(c-c_j)\shuffle d\|_{\ell_\infty,M_\epsilon}+\|c_j\shuffle (d-d_j)\|_{\ell_\infty,M_\epsilon}\\
&\leq K_\epsilon \|(c-c_j)\|_{\ell_\infty,M} \|d\|_{\ell_\infty,M}+K_\epsilon \|c_j\|_{\ell_\infty,M} \|(d-d_j)\|_{\ell_\infty,M}.
\end{align*}
Thus, $\lim_{j\rightarrow\infty}\|(c\shuffle d)-(c_j\shuffle d_j)\|_{{\ell_\infty,M_\epsilon}}=0$, proving the second part of the theorem.
\endpr

The next lemma will be needed for proving continuity of the shuffle inverse as well as for proving continuity of the composition product in the next section.

\begin{lemma} \label{le:shuffle-power-difference}
If $c$ and $c_j$, $j\geq 1$ are proper series in $\ell_{\infty,M}(X^\ast,\K)$ for some $M\in\nat$, and $\|c-c_j\|_{\ell_\infty,M}\rightarrow 0$ as $j\rightarrow \infty$,
then for $N\in\nat$ sufficiently large it follows that
\begin{equation*}
\sum_{n=1}^\infty \|c^{\shuffle n}-c_j^{\shuffle n}\|_{\ell_\infty,N}\rightarrow 0
\end{equation*}
as $j\rightarrow \infty$.
\end{lemma}

\begpr
 It is first shown that the sum is uniformly bounded for some large enough $N_0 \in \mathbb{N}$. The fact that $c_j \to c$ in $\ell_{\infty , M }$ and that $\|\cdot \|_{\ell_{\infty,N}} < \| \cdot \|_{\ell_{\infty,M}} $ for nonzero proper elements whenever $N > M$, implies that one can choose $N \in \mathbb{N}$ so that $ \| c \|_{\ell_{\infty,N}} + \sup_{j \in \mathbb{N}} \| c_j \|_{\ell_{\infty,N}} \leq \frac{1}{2} $. Define a proper series $d \in \ell_{\infty , N}(X^\ast,\K)$ by
\begin{equation*}
    (d , \eta) := (\| c \|_{\ell_{\infty,N}} + \sup_{j \in \mathbb{N}} \| c_j \|_{\ell_{\infty,N}} ) \, N^{| \eta |} \, | \eta |!, \; \; \eta \neq \emptyset.
\end{equation*}
Then $(d , \eta ) \geq | (c , \eta ) | + | (c_j , \eta) | $ for any $\eta \in X^*$ and all $j \in \mathbb{N}$. In fact, $(d^{\shuffle n} , \eta ) \geq | (c^{\shuffle n} , \eta ) | + | (c_{j}^{\shuffle n} , \eta ) |$ for any $n \geq 1$ and all $j \in \mathbb{N}$ by a standard induction argument. In particular, if $N' \geq N$, then
\begin{equation*}
    \| d^{\shuffle n} \|_{\ell_{\infty,N'}} \geq \| c^{\shuffle n} - c_{j}^{\shuffle n} \|_{\ell_{\infty,N'}}.
\end{equation*}
Now observe that
\begin{align*}
    (d^{\shuffle \, n}, \eta ) &\leq ( \| c \|_{\ell_{\infty , N}} + \sup_{j \in \mathbb{N}} \, \| c_j \|_{\ell_{\infty, N}} )^{n} \, N^{| \eta |} \, \binom{(n-1) + | \eta |}{n-1} \, | \eta |! \\
    &\leq \frac{1}{2^n} \, N^{| \eta |} \, \binom{(n-1) + | \eta |}{n-1} \, | \eta |! \\
    &= \frac{1}{2^n} \, (4N)^{| \eta |} \, | \eta |! \, \binom{(n-1) + | \eta |}{n-1} \, \frac{1}{4^{| \eta |}}.
\end{align*}
The first inequality can be shown via induction. Moreover, one can show the existence of a positive constant $K$ for which
\begin{equation*}
    \binom{(n-1) + | \eta |}{n-1} \, \frac{1}{4^{| \eta |}} \leq K, \; \; \forall n \in \mathbb N,\; |\eta| \geq n.
\end{equation*}
As $d$ is proper, $(d^{\shuffle n} , \eta ) = 0$ for all words $| \eta | < n$. Therefore,
\begin{equation*}
    (d^{\shuffle n} , \eta ) \leq \frac{1}{2^n} K (4N)^{| \eta |} \, | \eta |! , \; \; \forall \eta \in X^*,
\end{equation*}
so that $\| d^{\shuffle n} \|_{\ell_{\infty,4N}} \leq K/2^n$. Setting $N_0 := 4N$ gives
\begin{align*}
    \sum_{n=1}^\infty \|c^{\shuffle n}-c_j^{\shuffle n}\|_{\ell_\infty,N_{0}} \leq \sum_{n=1}^\infty \|d^{\shuffle n}\|_{\ell_{\infty,N_{0}}} \leq \sum_{n=1}^\infty \frac{1}{2^{n}} \, K = K < \infty.
\end{align*}
Having shown that the sum is uniformly bounded, it is now claimed that for each $n \geq 1$, $\lim_{j \to \infty} \| c^{\shuffle n} - c_{j}^{\shuffle n} \|_{4N} = 0$. If this holds, then
\begin{equation*}
    \lim_{j \to \infty} \sum_{n=1}^\infty \|c^{\shuffle n}-c_j^{\shuffle n}\|_{\ell_\infty,4N} = \sum_{n=1}^\infty  \lim_{j \to \infty} \|c^{\shuffle n}-c_j^{\shuffle n}\|_{\ell_\infty,4N} = 0.
\end{equation*}
To prove the claim, define $N_n := N (1 + (1 - \frac{1}{n} ) )$ for $n \geq 1$ so that $2N > N_n > N_{n-1} > \dots > N_1 = N$. It is shown by induction on $n \geq 1$ that $\| c^{\shuffle n} - c_{j}^{\shuffle n} \|_{\ell_{\infty,N_n}} \rightarrow 0$ as $j\to\infty$. The case $n = 1$ follows immediately as $N_1 = N$ and $\| c - c_j \|_{\ell_{\infty,N}} \leq \| c - c_j \|_{\ell_{\infty,M}}$. Let $n > 1$. Using the bilinearity of the shuffle product it follows that
\begin{align*}
    \| c^{\shuffle \, n} - c_{j}^{\shuffle \, n} \|_{\ell_{\infty,N_n}} &= \| (c - c_j) \shuffle c^{\shuffle \, (n-1)} + c_j \shuffle (c^{\shuffle \, (n-1)} - c_{j}^{\shuffle \, (n-1)}) \|_{\ell_{\infty,N_n}} \\
    &\leq K_n \, \| c - c_j \|_{\ell_{\infty,N_{n-1}}} \, \| c^{\shuffle \, (n-1)} \|_{\ell_{\infty,N_{n-1}}} \\
    &\hspace*{0.2in}+ K_n \, \| c_j \|_{\ell_{\infty,N_{n-1}}} \, \| c^{\shuffle \, (n-1)} - c_{j}^{\shuffle \, (n-1)} \|_{\ell_{\infty,N_{n-1}}},
\end{align*}
where as in Lemma~\ref{le:shuffle-Ke} the $K_n > 0$ are the constants corresponding to the $\epsilon_n > 0$ for which $N_n = (1 + \epsilon_n) N_{n-1}$. By the induction hypothesis, the latter expression tends to zero, implying the same for the former. This proves the claim since $\| c^{\shuffle n} - c_{j}^{\shuffle n} \|_{\ell_{\infty,4N}} \leq \| c^{\shuffle n} - c_{j}^{\shuffle n} \|_{\ell_{\infty,N_n}} , \; \; \forall n \geq 1 $.
\endpr

\begin{propose} \label{prop: ContShuffleInverse}
Denote by $(\allseriesLC)^\times$ the set of invertible elements of the algebra $(\allseriesLC, \shuffle)$. The shuffle inverse
\begin{equation*}
    \shuffle^{-1} \colon (\allseriesLC) ^\times \to (\allseriesLC) ^\times \; , \; c \mapsto (c , \emptyset)^{-1} \; \sum_{k \geq 0} (c^\prime)^{\shuffle k},
\end{equation*}
where $c' = \mathbf{1} - c/(c, \emptyset)$, is well defined and continuous.
\end{propose}

\begpr
Well definedness follows from \cite[Theorem 5]{Gray-etal_AUTO14}. To show continuity, first observe that $(\allseriesLC) ^\times$ is an open subset of $\allseriesLC$. Indeed it is easily verified that for any $\eta \in X^*$ the evaluation map $a_{\eta}$ is continuous on the Silva space $\allseriesLC$. In particular, $(\allseriesLC)^\times = a_{\emptyset} ^{-1} (\mathbb K \setminus \{0\} )$ is open. Since $\allseriesLC$ is sequential, the same is true for the open subset $(\allseriesLC) ^\times$, and consequently it suffices to test continuity of $\shuffle ^{-1}$ via sequences. With this in mind suppose $c_j \to c$ for elements $c_j , c \in (\allseriesLC) ^\times$, say $\| c - c_j \|_{\ell_{\infty, M}} \to 0$ for some $M > 0$. Then also $\| c' - c_{j} ' \|_{\ell_{\infty, M}} \to 0$. Since $c'$ and $c_j '$ are proper series, applying Lemma \ref{le:shuffle-power-difference} gives
\begin{equation*}
    \left\| \sum_{k =0}^\infty (c^\prime)^{\shuffle k} - (c_{j}^\prime)^{\shuffle k} \right\|_{\ell_{\infty, N}} =
    \left\| \sum_{k = 1}^\infty (c^\prime)^{\shuffle k} - (c_{j}^\prime)^{\shuffle k} \right\|_{\ell_{\infty, N}} \to 0
\end{equation*}
for some $N > M$. Hence,
\begin{equation*}
    \| \shuffle{^{-1}} (c) - \shuffle^{-1} (c_j) \|_{\ell_{\infty, N}} =
    \left\| (c, \emptyset)^{-1} \; \sum_{k =0}^\infty (c^\prime)^{\shuffle k} - \; (c_j , \emptyset)^{-1} \; \sum_{k=0}^\infty (c_{j}^\prime)^{\shuffle k} \right\|_{\ell_{\infty, N}} \to 0,
\end{equation*}
or in other words, $\shuffle^{-1}(c_j ) \to \shuffle^{-1} (c)$.
\endpr

\subsection{Continuity of the composition product}

In addition to Lemma \ref{le:shuffle-power-difference} the next result is needed in order to address the continuity of the composition product.

\begin{lemma} \label{le:composition-Ke}
Fix $M>0$. If $c\in\ell_{\infty,M}(X^\ast,\K^\ell)$ and $d\in\ell_{\infty,M}(X^\ast,\K^m)$, then $c\circ d\in \ell_{\infty,M_\epsilon}(X^\ast,\K^\ell)$
for any $M_\epsilon=M(1+\epsilon)$, $\epsilon>\phi(m \|d\|_{\ell_\infty,M})$ with $\phi(x)=x/2+\sqrt{x^2/4+x}$
and
\begdis
\|c\circ d\|_{\ell_\infty,M_\epsilon}\leq \|c\|_{\ell_\infty,M} (K_\epsilon\circ\phi)(m\|d\|_{\ell_\infty,M}),
\enddis
where
$K_\epsilon(a)=\sup_{\eta\in X^\ast}({\abs{\eta}+1})(1+a)^{\abs{\eta}}/(1+\epsilon)^{\abs{\eta}}$.
\end{lemma}

\begpr
It was shown in \cite{Gray-Li_05} that under the stated conditions
\begdi
 |(c\circ d,\eta )| \leq  \|c\|_{\ell_\infty,M} ((1+\phi(m \|d\|_{\ell_\infty,M})) M)^{|\eta|}(|\eta|+1)!,\;\; \forall \eta \in X^\ast.
\enddi
Therefore,
\begdi
\frac{|(c\circ d,\eta)|}{{M_\epsilon}\abs{\eta}!} \leq  \|c\|_{\ell_\infty,M} (1+\phi(m \|d\|_{\ell_\infty,M}))^{|\eta|}\frac{(|\eta|+1)}{(1+\epsilon)^{\abs{\eta}}},\;\; \forall \eta \in X^\ast.
\enddi
Taking the supremum over $X^\ast$ proves the lemma.
\endpr

\begin{theorem} \label{th:Silva-continuity-comp-product}
The composition product on $\allseriesmLC$ is continuous in the Silva topology.
\end{theorem}

\begpr
Left and right continuity of the composition product is first proved,
beginning with left continuity. Let $M>0$ be fixed. Let $c,d\in\ell_{\infty,M}(X^\ast,\K^m)$,
and assume $c_j$, $j\geq 1$ is a sequence in $\ell_{\infty,M}(X^\ast,\K^m)$ converging to $c$.
Applying Lemma~\ref{le:composition-Ke} gives
\begin{align*}
\|(c\circ d)-(c_j\circ d)\|_{{\ell_\infty,M_\epsilon}}&=\|(c-c_j)\circ d\|_{{\ell_\infty,M_\epsilon}} \\
&\leq  \|c-c_j\|_{\ell_\infty,M} (K_\epsilon\circ\phi)(m\|d\|_{\ell_\infty,M}).
\end{align*}
Thus, $\lim_{j\rightarrow \infty} \|(c\circ d)-(c_j\circ d)\|_{{\ell_\infty,M_\epsilon}}=0$.

Right continuity is addressed next. It
is more complicated given the nonlinearity in the right argument of the product.
Let $c,d\in\ell_{\infty,M}(X^\ast,\K^m)$
and assume $d_j$, $j\geq 1$ is a sequence in $\ell_{\infty,M}(X^\ast,\K^m)$ converging to $d$.
For a fixed $\xi\in X^\ast$, observe that
\begin{align*}
|((c\circ d)-(c\circ d_j),\xi)|&=\left|\sum_{\eta\in X^\ast} (c,\eta)(\eta\circ d-\eta\circ d_j,\xi)\right| \\
&\leq \sum_{n=0}^\infty \|c\|_{\ell_\infty,M} M^n\,n! \left| \sum_{\eta\in X^n} (\eta\circ d-\eta\circ d_j,\xi)\right| \\
&= \|c\|_{\ell_\infty,M} \sum_{n=0}^\infty M^n\,n! \left| \sum_{r_0\geq 0,\ldots,r_m\geq 0 \atop r_0+\cdots+r_m=n} \right.
((x_0^{r_0}\shuffle\cdots\shuffle x_m^{r_m})\circ d \\
&\hspace*{0.2in}\left.\rule{0in}{0.4in} -(x_0^{r_0}\shuffle\cdots\shuffle x_m^{r_m})\circ d_j,\xi)\right|. \\
\end{align*}
Applying the identities $x_i^{\shuffle n}=n!\,x_i^n$, $n\geq 0$ and $(c\shuffle d)\circ e=(c\circ e)\shuffle (d\circ e)$ gives
\begin{align*}
|((c\circ d)-(c\circ d_j),\xi)|
&\leq \|c\|_{\ell_\infty,M} \sum_{n=0}^\infty M^n \left| \sum_{r_0\geq 0,\ldots,r_m\geq 0 \atop r_0+\cdots+r_m=n}
{n \choose r_0\cdots r_m}
((x_0^{\shuffle r_0}\shuffle\cdots\shuffle x_m^{\shuffle r_m})\circ d \right. \\
&\hspace{0.2in} \left. \rule{0in}{0.4in} -(x_0^{\shuffle r_0}\shuffle\cdots\shuffle x_m^{\shuffle r_m})\circ d_j,\xi)\right|. \\
&=\|c\|_{\ell_\infty,M} \sum_{n=0}^\infty M^n \left|\left(\left(\sum_{k=0}^m x_k\circ d\right)^{\shuffle n}-\left(\sum_{k=0}^m x_k\circ d_j\right)^{\shuffle n},\xi\right) \right| \\
&=\|c\|_{\ell_\infty,M} \sum_{n=0}^\infty \left|\left(\bar{d}^{\shuffle n}-\bar{d}_j^{\shuffle n},\xi\right) \right|,
\end{align*}
where $\bar{d}:=Mx_0\sum_{k=0}^m d[k]$ and $\bar{d_j}:=Mx_0\sum_{k=0}^m d_j[k]$ are proper series in $\allseries$.
Here $d[k]$ denotes the $k$-th component series of $d$. It is clear that $\bar{d}\in \ell_{\infty,M}(X^\ast,\K)$, and
$\bar{d}_j$ is a sequence in $\ell_{\infty,M}(X^\ast,\K)$. Furthermore,
$\bar{d}_j\rightarrow \bar{d}$ as $j\rightarrow \infty$ since
\begin{align*}
\|\bar{d}-\bar{d}_j\|_{\infty,M}&=\sup_{\eta\in X^\ast} \frac{|(\bar{d}-\bar{d}_j,\eta)|}{M^{\abs{\eta}}\abs{\eta}!} \\
&\leq \sum_{k=1}^m \sup_{\eta\in X^\ast} \frac{M|(x_0(d[l]-d_j[k]),\eta)|}{M^{\abs{\eta}}\abs{\eta}!} \\
&= \sum_{k=1}^m \sup_{x_0\eta\in X^\ast} \frac{|(d[k]-d_j[k],\eta)|}{M^{\abs{\eta}}(\abs{\eta}+1)!} \\
&= \sum_{k=1}^m \sup_{\eta\in X^\ast} \frac{|(d[k]-d_j[k],\eta)|}{M^{\abs{\eta}}\abs{\eta}!(\abs{\eta}+1)} \\
&\leq m \|d-d_j\|_{\ell_\infty,M}.
\end{align*}
Finally, right continuity follows by applying Lemma~\ref{le:shuffle-power-difference} with $N>M$ sufficiently large so
that
\begin{align}\label{eq:comp:rcont}
\|(c\circ d)-(c\circ d_j)\|_{\ell_\infty,N}&\leq \|c\|_{\ell_\infty,M}
\sum_{n=0}^\infty \|\bar{d}^{\shuffle n}-\bar{d}_j^{\shuffle n}\|_{\ell_\infty,N}.
\end{align}

Note that the estimates for left and right continuity imply joint continuity of the composition product due to the following simple observation that for $N$ as above
\begin{align*}
 \|(c\circ d)-(c_j\circ d_j)\|_{{\ell_\infty,N}} &\leq \|(c\circ d)-(c_j \circ d)\|_{{\ell_\infty,N}} + \| (c_j \circ d) -(c_j\circ d_j)\|_{{\ell_\infty,N}}\\
 & \stackrel{\eqref{eq:comp:rcont}}{\leq}  \|(c-c_j)\circ d\|_{{\ell_\infty,N}} +  \|c_j\|_{\ell_\infty,M}
\sum_{n=0}^\infty \|\bar{d}^{\shuffle n}-\bar{d}_j^{\shuffle n}\|_{\ell_\infty,N},
\end{align*}
where the last inequality is a direct consequence of \eqref{eq:comp:rcont}. Hence we see that the product is sequentially continuous (as each sequence convergent in the Silva topology is already contained in one of the Banach steps). By Lemma \ref{lem:Silvaprop} implies that the product is continuous as each of the $p_{M,K}$ is continuous for every $M,K >0$.
\endpr

\section{Analyticity of the composition and shuffle product}

In this section it is proved that the formal power series products and inverse presented in the previous sections are not only continuous but also analytic.
Note that on the infinite-dimensional spaces involved, both complex and real analyticity make sense, cf.\ Appendix \ref{app:smooth}.
For real analyticity one needs only to identify the complexification of the spaces $\allseriesR$ and $\allseriesLCR$.

As locally convex spaces, the complexification of $\allseriesR$ is $\allseriesC$. This is clear on the level of vector spaces,
and for the topology simply note that as topological vector spaces $\allseriesC = \allseriesR \oplus i\allseriesR$.
Similarly, the complexification of the Silva space $\allseriesLCR$ is $\allseriesLCC$. Again this is clear on the level of vector spaces but more complicated on the level of the vector space topologies. However, also the vector space topologies coincide as it is easy to see that for every $M >0$ the Banach space $\ellInfty{M}{X^*}{\C^\ell}$ is the complexification of $\ellInfty{M}{X^*}{\R^\ell}$, and the inductive limit of a sequence of compact operators between Banach spaces commutes with the formation of complexifications \cite[Theorem 3.4]{MR1878717}.

Having identified the complexification of the infinite-dimensional spaces, observe that the shuffle product, the composition product and the shuffle inverse are all well defined on both the complexification and on the real space. Hence, if it can be proved that these mappings are holomorphic on the complexification, then real analyticity is obtained for the corresponding mappings on the real space. Before continuing with the shuffle product and the shuffle inverse, it is helpful to recall a special type of locally convex algebra.

\begin{definition}
 Let $(A,\beta)$ be an associative unital locally convex algebra, i.e.,\ $A$ is a locally convex space such that the bilinear map $\beta$ is continuous and admits a unit $\mathbf{1}$ with $\beta(\mathbf{1},x) =x = \beta(x,\mathbf{1})$. Then $A$ is called a \bfem{continuous inverse algebra} (CIA) if the unit group $A^\times$ is an open subset of $A$, and inversion $\iota \colon A^\times \rightarrow A^{\times}$ is continuous.
\end{definition}

\begin{propose}
 The algebras  $(\allseries, \shuffle)$ and $(\allseriesLC, \shuffle)$ are continuous inverse algebras with respect to their natural topologies.
\end{propose}

\begpr
It was shown in Theorem \ref{thm:shuffle:cont} that the bilinear shuffle product is continuous with respect to the Silva and the \Frechet topology. Furthermore, the non proper series are precisely the invertible elements with respect to the shuffle product. By definition of a non proper series it is evident that if $A$ is either the algebra $\allseries$ or the algebra $\allseriesLC$, then $A^\times = a_{\emptyset}^{-1} (\C \setminus \{0\})$ is open as the preimage of an open set under a continuous map.
Continuity of the shuffle inverse for the Silva topology on $\allseriesLC$ was established in Proposition \ref{prop: ContShuffleInverse}. To see that the shuffle inverse is also continuous on $(\allseries)^\times$ it suffices to test continuity of the composition $a_\eta \circ \shuffle^{-1}$ for every $\eta \in X^*$. However, due to the definition of the shuffle inverse, it is clear that $a_\eta \circ \shuffle^{-1}(c)$ is a polynomial in finitely many evaluations of the series $c$.  Therefore, $a_\eta \circ \shuffle^{-1}$ is continuous, and hence the shuffle inverse is continuous in the \Frechet topology.
\endpr

It is well known that the unit group of a CIA is an infinite-dimensional Lie group. Before stating the next result, recall the following notion from infinite-dimensional Lie theory.

\begin{definition}
 Consider a Lie group $G$ with unit $1$ and write $\Lf(G)$ for the Lie algebra of $G$. Let $\lambda_g \colon G\rightarrow G, \lambda_g(h)=gh$ be the left multiplication with a fixed element $g \in G$. Then $G$ is called \emph{$C^r$-regular}, $r\in \mathbb{N}_0\cup\{\infty\}$, if for each $C^r$-curve $u\colon [0,1]\rightarrow \Lf(G)$ the initial value problem
 \begin{displaymath}
  \begin{cases}
  \dot \gamma(t)= \gamma(t) . u(t) \coloneq T\lambda_{\gamma(t)}(u(t))\\ \gamma(0) = 1
  \end{cases}
 \end{displaymath}
 has a (necessarily unique) $C^{r+1}$-solution
 $\Evol (u)\coloneq\gamma\colon [0,1]\rightarrow G$ and the map
 \begin{displaymath}
  \evol \colon C^r([0,1],\Lf(G))\rightarrow G,\quad u\mapsto \Evol (u)(1)
 \end{displaymath}
 is smooth.\footnote{The function space $C^r([0,1],\Lf(G))$ is endowed with the compact open $C^r$-topology (controlling a function and its derivatives on compact subsets). With this topology and pointwise addition and scalar multiplication $C^r([0,1],\Lf(G))$ is a locally convex space. Thus, it makes sense to define smooth mappings on this space, cf.\, Appendix \ref{app:smooth}.} A $C^\infty$-regular Lie group $G$ is called \emph{regular}
 \emph{(in the sense of Milnor}).
\end{definition}
Every Banach Lie group is $C^0$-regular (cf. \cite{neeb2006}). Several important results in infinite-dimensional Lie theory are only available for regular Lie groups. For example the interplay between Lie algebra and Lie group hinges on regularity as this property guarantees existence of a smooth Lie group exponential function. Moreover, if one wants to lift morphisms of Lie algebras to the Lie group by integration, this requires the group to be regular, cf.\ \cite{KM97}.
\begin{propose}
 The group $((\allseries)^\times, \shuffle)$ with the \Frechet topology and the group\\ $((\allseriesLC)^\times,\shuffle)$ with the Silva topology are $C^0$-regular analytic Lie groups.
\end{propose}

\begpr
It was established that the groups are unit groups of continuous inverse algebras, hence they are infinite-dimensional analytic Lie groups by \cite[Theorem 5.6]{Glo02}. Moreover, since the shuffle product is abelian, and $\allseries$ and $\allseriesLC$ are both complete locally convex spaces, an application of \cite[p.3 Corollary and Proposition 3.4 (a)]{GN12} shows that the Lie groups $\allseries$ and $\allseriesLC$ are $C^0$-regular (even with analytic evolution map $\evol$).
\endpr

\begin{rem}
In \cite[Lemma 2.2]{GN12} it was proved that the solution to the initial value problem for regularity in the unit group of a CIA is given by the Volterra series
\begin{equation}
\gamma(t) = 1 + \sum_{n=1}^\infty \int_0^t \int_0^{t_{n-1}} \cdots \int_0^{t_2} \eta (t_1) \cdots \eta (t_n) \mathrm{d}t_1 \ldots \mathrm{d}t_n. \label{eq:Volterra}
\end{equation}
Hence the Volterra series describes both the solution of the initial value problem in $\allseries$ and the subgroup $\allseriesdeltaLC$.
\end{rem}

\begin{propose}\label{prop:compose:analytic}
 The composition product on $\allseries$ and on $\allseriesLC$ is analytic.
\end{propose}

\begpr
In light of the previous observations regarding the complexifications, it suffices to prove the statement for the case $\K = \C$.
 Fix $\eta \in X^*$. By definition of the composition product, an induction argument shows that $a_\eta (c\circ d)$ is a polynomial in finitely many $a_\gamma(c)$ and $a_\rho (d)$ for words such that $|\gamma|,|\rho| \leq |\eta|$ (for a detailed proof see \cite[Lemma 83]{palm20}). As the coordinate functions are continuous linear (thus holomorphic) in the \Frechet topology on $\allseries$ and in the Silva topology on $\allseriesLC$, one can deduce the following:
\begin{enumerate}
 \item the composition product $\circ \colon \allseries^2 \rightarrow \allseries$ is continuous with respect to the \Frechet topology (which is initial with respect to the $a_\eta$);
 \item for every $\eta \in X^*$ the map $(c,d) \mapsto a_\eta (c\circ d)$ is holomorphic both on $\allseries^2$ and on $\allseriesLC^2$.
\end{enumerate}
Furthermore, the coordinate functions $a_\eta, \eta \in X^*$ separate the points on $\C\langle\langle X \rangle\rangle$ and on $\C_{LC}\langle\langle X \rangle\rangle$.
Now apply Lemma \ref{lem:folklore}. Since the composition product is continuous on $\C\langle\langle X \rangle\rangle$ and analytic after composition with $a_\eta$, $\eta \in X^*$, the composition product is analytic as a mapping on $\C\langle\langle X \rangle\rangle$. A similar argument holds for the composition product on $\C_{LC}\langle\langle X \rangle\rangle$ as continuity for this product was established in Theorem \ref{th:Silva-continuity-comp-product}.
\endpr

\section{The Lie group $(\delta + \allseriesmLC,\circ,\delta)$}

A Lie group structure on the group $(\delta + \allseriesmLC,\circ ,\delta)$ is presented in this section.
This group is known to have an associated graded and connected Hopf algebra $(H,\mu,\Delta)$ as
described in Section~\ref{subsec:system-interconnections} and more completely in \cite[Section 3]{Gray-etal_SCL14}. This structure will play an important role in the proof of the Lie group property. Note, however, that $(\allseriesdeltam,\circ ,\delta)$ is \emph{not} the character group of said Hopf algebra, and thus the Lie theory for such groups from \cite{BDS16,DS18} is not directly applicable. The main claim, as stated below, is established from first principles.

\begin{theorem} \label{thm:Lieloc}
The group $(\delta + \allseriesmLC,\circ,\delta)$ is an analytic Lie group under the Silva topology.
\end{theorem}

\begpr
The proof is carried out in four main steps.

\textbf{Step 1:} \emph{The group product is continuous in the Silva topology}.
Fix $M\geq 0$ and let $M_\epsilon=M(1+\epsilon)$. If $c,d\in \ell_{\infty,M}(X^\ast,\K^m)$ then the proof
of Theorem~\ref{th:Silva-continuity-comp-product} can be easily modified to show that $c\modcomp d_\delta$ is
continuous in the Silva topology. Specifically, the only change is in the definition of $\bar{d}$ and $\bar{d}_j$.
For example, $\bar{d}=M(\sum_{k=0}^m x_k+x_0d[k])$. In which case, it follows directly that
$c_\delta\circ d_\delta=\delta+d+c\modcomp d_\delta$ is continuous in both its left and right arguments in
the Banach space $\ell_{\infty,M_\epsilon}(X^\ast,\K^m)$. Joint continuity follows then verbatim as in the proof of Theorem~\ref{th:Silva-continuity-comp-product}.\\[.5em]

\textbf{Step 2:} \emph{The group inverse is degreewise a polynomial}. Assume without loss of generality that $m=1$.
Let $c_j\rightarrow c$ in $\ellInfty{M}{X^*}{\K}$.
It was shown in \cite{Gray-etal_SCL14} that the composition inverse preserves local convergence. Thus, there exists an $M_1>0$ such that
$c_\delta^{\circ -1}\in\delta+\ellInfty{M_1}{X^*}{\K}$ and $(c_{\delta,j})^{\circ -1}\in\delta+\ellInfty{M_1}{X^*}{\K}$ for every $j\geq 1$.
Set $M_2=\max(M,M_1)$. Since $H$ is graded and connected with respect to the degree grading, it follows from Lemma~\ref{le:antipode-is-group-inverse} (cf.\, \cite{MR2523455}) that
\begin{align}
 (c_\delta^{\circ -1},\eta) =& S(a_\eta)(c) = -a_\eta(c)-\sum S(a_{(\eta_1)}^\prime)(c)a_{(\eta_2)}^\prime(c) \notag\\ =& -a_\eta(c)+\sum_{k=1}^{\deg(a_\eta)} (-1)^{k+1}\mu_k \circ \Delta_k^\prime (a_\eta)(c), \label{inverse:poly}
\end{align}
where  $\Delta^\prime a=\Delta a-a\otimes \mbf{1}_\delta-\mbf{1}_\delta\otimes a=\sum a^\prime_{(\eta_1)}\otimes a^\prime_{(\eta_2)}$ is the reduced
coproduct in the notation of Sweedler\footnote{Given the bijection between $\delta+\allseries$ and $\allseries$, for brevity $a_\eta(c_\delta)$ will be written as $a_\eta(c)$.}, $\Delta_k^\prime = \Delta_{k-1}^\prime \otimes \id$ is defined inductively, and $\mu_k$ is the $k$-fold multiplication in the target algebra. In particular, $a^\prime_{(\eta_1)}\in V_{n_1}$ and $a^\prime_{(\eta_2)}\in H_{n_2}$ with $n_1,n_2<n$. As the summation in \eqref{inverse:poly} is always finite, the $\eta$ component of $c^{\circ -1}_\delta$ is a polynomial in the variables $\{a_\xi(c):\deg(a_\xi)\leq \deg(a_\eta)\}$. This implies immediately that inversion is continuous (and analytic) in the \Frechet space $\delta + \allseriesmLC$. However, this does not yet yield continuity with respect to the Silva space topology on $\delta + \allseriesmLC$.
\\[.5em]

\textbf{Step 3:} \emph{Continuity of the group inverse in the Silva topology}.
It is first proved that inversion is continuous at the unit $\delta$. It is again assumed without loss of generality
that $m=1$. Recalling that $c_\delta:=\delta+c$, the series
$c_{\delta,j}= \delta + c_{j}, j \in \mathbb{N}$ converges to $\delta$ in the Silva topology if and only if the series
$c_j$ converges to $0$ in $\ellInfty{M}{X^*}{\K}$ for some $M>0$.
Fix $c\in\ellInfty{M}{X^*}{\K}$ and define $\bar{c}=\sum_{\eta\in X^\ast}KM^{\abs{\eta}}\abs{\eta}!\,\eta$ with $K= \|c\|_{{\ell_\infty,M}}$
so that $\abs{(c,\eta)}\leq (\bar{c},\eta)$, $\forall \eta\in X^\ast$.
It can be verified directly that $y=F_{\bar{c}_\delta}[u]=u+F_{\bar{c}}[u]$ has the state space realization
\begdi
\dot{z}=\frac{M}{K}(1+u),\;\;z(0)=K,\;\;y=z+u.
\enddi
Therefore, $y=F_{\bar{c}^{\circ -1}_\delta}[u]=u+F_{\bar{c}^{\circ -1}}[u]$ has the realization
\begeq \label{eq:inverse-realization}
\dot{z}=\frac{M}{K}(z^2-z^3)+z^2u,\;\;z(0)=K,\;\;y=-z+u.
\endeq
It is shown in \cite[Theorem 6]{Gray-etal_SCL14} that $c^{\circ -1}=(-c)@\delta$, where the right-hand side denotes
the generating series for the unity feedback system $v\mapsto y$ defined by $y=F_{-c}[u]$ and $u=v+y$. Combining this fact with a minor extension of \cite[Lemma 10]{Thitsa-Gray_SIAM12},
it follows that the condition $\abs{(c,\eta)}\leq (\bar{c},\eta)$ implies $\abs{(c^{\circ -1},\eta)} \leq \abs{(\bar{c}^{\circ -1},\eta)}$, $\forall \eta\in X^\ast$.
The fastest growing
coefficients of $\bar{c}^{\circ -1}$ have been shown to be the sequence $(\bar{c}^{\circ -1},x_0^k)$, $k\geq 0$ \cite[Lemma 7]{Thitsa-Gray_SIAM12}.
Therefore, for any word $\eta\in X^\ast$ of length $k$
\begdi
\abs{(c^{\circ -1},\eta)} \leq \abs{(\bar{c}^{\circ -1},\eta)}
\leq\abs{(\bar{c}^{\circ -1},x_0^k)}
=\abs{L_{g_0}^k h(z_0)},
\enddi
where the right-most inequality follows from \rref{eq:c-from-Lgh}
with $g_0(z)=(M/K)(z^2-z^3)$, $h(z)=-z$, and $z_0=K$ as derived in
\rref{eq:inverse-realization}. A direct calculation gives
\begeq \label{eq:first-barc-inverse-bound}
(\bar{c}^{\circ -1},x_0^k)=b_k(K)KM^k k!, \;\;k\geq 0,
\endeq
where the first few polynomials $b_k(K)$ are:
\begin{align*}
b_0(K)&=-1 \\
b_1(K)&=-1 +K  \\
b_2(K)&=-2  +5 K -3 K^2  \\
b_3(K)&=-6  +26 K -35 K^2 +15 K^3  \\
b_4(K)&=-24  +154 K -340 K^2 +315 K^3 -105 K^4  \\
b_5(K)&=-120  +1044 K -3304 K^2 +4900 K^3 -3465 K^4 +945 K^5  \\
b_6(K)&=-720  +8028 K -33740 K^2 +70532 K^3 -78750 K^4 +45045 K^5 -10395 K^6  \\
b_7(K)&=-5040 +69264 K -367884 K^2 +1008980 K^3 -1571570 K^4 +1406790 K^5 \\
&\hspace*{0.25in}-675675 K^6 +135135 K^7 \\
&\hspace*{0.08in}\vdots
\end{align*}
When $K\leq 1$ it is known that $b_k(K)\leq \bar{b}_k$, where $\bar{b}_k$, $k\geq 0$ is the integer sequence A112487 in \cite{OEIS},
namely, 1, 2, 10, 82, 938, 13778, 247210, \ldots.
Its exponential generating function is the real analytic function
\begdi
G(x)=\frac{-1}{1+W(-2\exp(x-2))},
\enddi
where $W$ is the Lambert W-function (see \cite[Example 5]{Thitsa-Gray_SIAM12}). In which case, there exists growth constants $\bar{K},\bar{M}>0$
such that $\bar{b}_k\leq \bar{K}\bar{M}^k k!$, $k\geq 0$. Combining this inequality with \rref{eq:first-barc-inverse-bound} gives
\begdi
\abs{(c^{\circ -1},\eta)}\leq  \|c\|_{{\ell_\infty,M}}\bar{K} (M\bar{M})^{\abs{\eta}} \abs{\eta}!,\;\;\forall \eta\in X^\ast.
\enddi
Hence, if $c_{\delta , j} \rightarrow \delta$ in $\K \times \ellInfty{M}{X^*}{\K}$, then $c_{\delta,j}^{\circ -1} \rightarrow \delta$ in $\ellInfty{M\overline{M}}{X^*}{\K}$.
Therefore, inversion is continuous at the unit with respect to the Silva topology.
Exploiting the fact that inversion is a group antimorphism, this implies that inversion is continuous
everywhere on $\delta + \allseriesmLC$ in the
Silva topology.\footnote{Alternatively, continuity can be deduced from a
more general criterion, see \cite[Lemma 1.3]{AaR05}.}\\[.5em]

\textbf{Step 4:} \emph{Group product and inverse are analytic.}
Since the complexification of $\delta+\re^m_{LC}\langle\langle X \rangle\rangle$ is $\delta+\C^m_{LC}\langle\langle X \rangle\rangle$, it suffices to consider the complex case. In view of Lemma \ref{lem:folklore} and Step 1, all one needs to prove is that for every $\eta \in X^\ast$ the mappings $(c_\delta,d_\delta) \mapsto a_\eta (c_\delta\circ d_\delta)$ and $c_\delta \mapsto a_\eta (c_\delta^{\circ -1})$ are holomorphic.
Regarding the composition product recall that $(\delta + c )\circ (\delta +d) = \delta +d + c\tilde{\circ} d_\delta$. Now for the mixed composition $\tilde{\circ}$ it was shown in the proof of Proposition \ref{prop:compose:analytic} that $a_\eta (c\tilde{\circ} d_\delta)$ is given by a polynomial in finitely many of the variables $a_\xi (c)$ and $a_\nu (d)$. Hence, this part of the product is analytic on $\delta+\C^m_{LC}\langle\langle X \rangle\rangle$, and therefore the composition product is analytic.
Similarly, for the inversion $\iota$, Step 2 shows that $a_\eta \circ \iota (c)$ is given as a polynomial in finitely many evaluations of $c$. As before, the coordinate functions are holomorphic and this implies that $a_\eta \circ \iota $ is holomorphic on $\delta+\C^m_{LC}\langle\langle X \rangle\rangle$. Hence, the inversion is also holomorphic.
\endpr

The argument for the Lie group structure on subsets of locally convergent series can be adapted almost verbatim to the case where no convergence of the series is assumed.

\begin{corollary}
 The group $(\delta + \allseriesm , \circ , \delta)$ is an analytic Lie group.
\end{corollary}

\begpr
Again it suffices to prove the case where $\K = \C$. In Step 2 of the proof for Theorem \ref{thm:Lieloc} it was shown that after composition with a coordinate function $a_\eta$ both the composition and the inversion in the group are given by a polynomial in finitely many coordinate functions applied to the arguments. Since the \Frechet
topology is initial with respect to the coordinate functions, it follows directly that the group operations are continuous. Applying Lemma \ref{lem:folklore} gives immediately that the group operations are also analytic.
\endpr

While the \Frechet Lie group $\delta + \allseriesm$ is much simpler (topologically speaking) than the Silva group $\delta+\allseriesmLC$, it supplies a useful template for the Lie theoretic arguments considered next, namely, identifying the Lie algebra and proving that a Lie group is regular in the sense of Milnor. The first goal is to establish these properties for the simpler \Frechet Lie group. Subsequently, it is shown that these results then imply corresponding properties for the Silva Lie group. However, it is first necessary to introduce a new structure which will yield a convenient description of the Lie bracket. This structure is the so called {\em pre-Lie product}, which was developed in \cite{Foissy_15} for the case where $m=1$ and generalized in \cite[Section 3.2]{DuffautEspinosa2016609} for the case where $m\geq 1$.

\begde  \label{de:Foissy-pre-Lie}
Let $X=\{x_0,x_1,\ldots,x_m\}$ and denote by $d[i]$ the $i$-th component of a series
$d\in\allseriesm$.  The pre-Lie product is
the bilinear product on $\allseriesm\times\allseriesm$
\begdi
c\lhd d=\sum_{\eta\in X^\ast} (c,\eta)\,\eta\lhd d,
\enddi
where $\eta\lhd d$ is defined inductively by
\begin{align*}
(x_0\eta)\lhd d &=x_0(\eta\lhd d) \\
(x_j\eta)\lhd d &=x_j(\eta\lhd d)+x_0(\eta\shuffle d[j]),\quad j=1,2,\ldots,m
\end{align*}
and $\emptyset\lhd d=0$.
\endde

This product can be viewed
as the linear part of the group product, that is,
\begeq \label{eq:group-product-expansion}
c_\delta\circ d_\delta=\delta+d+c\modcomp d_\delta=\delta+c+d+c\lhd d+{\mathcal O}(c,d^2),
\endeq
where ${\mathcal O}(c,d^2)$ denotes all terms depending linearly on $c$ and on higher powers of $d$.
One can show that the pre-Lie product preserves the length of words in the sense that
$(\eta\lhd\xi,\nu)=0$ when $\abs{\eta}+\abs{\xi}\neq \abs{\nu}$. Therefore, the product is well defined as it is locally finite.
Moreover, defining $d[0]=0$, the recursive formulas reduce to a single expression
\begin{align}\label{preLieformula}
 (x_j\eta)\lhd d &=x_j(\eta\lhd d)+x_0(\eta\shuffle d[j]), \quad j \in \{0,1,\ldots , m\}.
\end{align}

\begin{example}\label{ex:preLie}
{\rm
Consider the computation of the pre-Lie product for a few words of short length.
For example, if $c = x_0^n, n \in \mathbb{N}$ and $d \in \allseriesm$, then $x_0^n \lhd d = x_0^n (\emptyset \lhd d) = 0$.
For any $x_k\in X$ with $k\neq 0$,
\begin{align} \label{eq:lhd1}
x_k \lhd d = x_k (\emptyset \lhd d) +x_0 (\emptyset \shuffle d[k]) = x_0 d[k],\quad  k \in \{0, 1,\ldots ,m\}.
\end{align}
Observe $a_{x_k}(x_k \lhd d)=0$ as every word in the support of $x_k \lhd d$ must have the prefix $x_0$.
Furthermore, it is clear from \eqref{eq:lhd1} that the length of the words in $\supp(x_k \lhd d)$ coincide with the length of those
in $\supp(d)$ except incremented by one. On the other hand, if $d= \mathbf{e}_k \eta$ (where $\mathbf{e}_k \in \mathbb{R}^m$ is the $k$-th unit vector), then
$$\deg ( a_{x_k \lhd d}) = 2 + \deg (a_\eta) \geq \deg(a_{x_k})+ \deg (a_\eta), \quad k \in \{1,2,\ldots,m\}.$$
Indeed one always obtains $|\eta \lhd d|+|\eta\lhd d|_{x_0} \geq |\eta| +|\eta|_{x_0} + |d| +|d|_{x_0}$ (where the length of a sum of words is defined as the maximum of the lenght of the words).
Consider next $c=x_jx_k$ where both $j$ and $k$ are not zero. Applying the definition gives
\begin{align}
x_j x_k \lhd d &= x_j(x_k\lhd d) + x_0 (x_k \shuffle d[j]) \nonumber \\
&= x_j(x_k (\emptyset \lhd d) + x_0 d[k])+x_0 (x_k \shuffle d[j])\nonumber \\
&= x_jx_0d[k]+x_0(x_k \shuffle d[j]).\label{eq:lhd2}
\end{align}
For comparison, it follows from \rref{eq:mixed-composition-product} that
\begin{align*}
x_jx_k\modcomp d_\delta&=\phi_d(x_jx_k)(\mbf{1}) \\
&=\phi_d(x_j)\circ \phi(x_k)(\mbf{1}) \\
&=\phi_d(x_j)(x_k+x_0d[k]) \\
&=x_j(x_k+x_0d[k])+x_0(d[j]\shuffle(x_k+x_0d[k])) \\
&=x_jx_k+x_jx_0d[k]+x_0(d[j]\shuffle x_k)+x_0(d[j]\shuffle(x_0d[k])) \\
&=x_jx_k+x_jx_k\lhd d+x_0(d[j]\shuffle(x_0d[k])),
\end{align*}
which is consistent with \rref{eq:group-product-expansion}.
Applying now the coordinate function $a_{x_jx_k}$ to \eqref{eq:lhd2} gives $a_{x_jx_k}(x_jx_k \lhd d)=0$ for any series $d$. A trivial induction shows that
$$a_\eta (\eta \lhd d) = 0,\; \forall \eta \in X^\ast,\; d \in \allseriesm.$$
Finally, consider a word $\eta$ with $|\eta|_{x_0}=0$. Observe $a_\eta (\rho \lhd d)=0$ because every word in the support of $a_\eta (\rho \lhd d)$ must contain at least one $x_0$ and $\abs{\eta}_{x_0}=0$.%consider a word $\eta = \prod_{k=1}^r x_{j_k}$ such that all $j_k\neq 0$. Then for all $\rho = \prod_{i=1}^s x_{\ell_i}\in X^\ast$ with $|\rho|\leq |\eta|$ and $d \in \allseriesm$
%$$a_\eta (\rho \lhd d) = \delta_{j_1,\ell_1} a_{\prod_{k=2}^r x_{j_k}} \left(\prod_{i=2}^s x_{\ell_i} \lhd d\right) = \cdots = \left(\prod_{k=1}^{\min (r,s)} \delta_{j_k,\ell_k} \right) a_{\prod_{k=\min (r,s)}^n x_{j_k}}(\emptyset \lhd d) =0,$$
%where $\delta_{i,j}$ denotes the Kronecker delta function.
}
\end{example}

\begin{propose}
The Lie algebra of $\delta + \allseriesm$ is the space $\allseriesm$ with the Lie bracket given by the formula
\begin{align}
\LB[c,d] = c \lhd d-  d\lhd c \label{compositionLiebracket}.
\end{align}
\end{propose}

\begpr
The Lie bracket of the Lie algebra associated to the Lie group $\delta + \allseriesm$ is given by evaluating the Lie bracket of left invariant vector fields on $\delta + \allseriesm$ at the identity $\delta$.
Note that since $\delta + \allseriesm$ is an affine subspace of $\allseriesdeltam$, it is easy to see that the left-invariant vector field associated to $c \in \allseriesm$ is given by the formula
$X^c (\delta + e) = c + e \lhd c $, hence
\begin{align*}
\LB[c,d] &= \LB[X^c,X^d](\delta) \\
&= (dX^d\circ X^c - dX^c \circ X^d)(\delta) \\
&= \lim_{t \to 0} t^{-1} \, (X^d (\delta + tc) - X^c (\delta) - X^c (\delta + td) + X^d (\delta) )\\
&= \lim_{t \to 0} t^{-1} \, (d + tc \lhd d - d - c - td \lhd c + c ) \\
&= \lim_{t \to 0} t^{-1} \, t \, (c \lhd d - d \lhd c ) = c \lhd d - d \lhd c.
\end{align*}
\endpr

\begin{corollary}
 The Lie algebra of $(\delta + \allseriesmLC,\circ,\delta)$ is $\allseriesmLC$ with bracket \eqref{compositionLiebracket}.
\end{corollary}

\begpr
The canonical inclusion $\iota \colon \delta + \allseriesmLC \rightarrow \delta + \allseries$ is the restriction of a continuous linear map to a closed (affine linear) subset, whence smooth. Obviously it is a Lie group morphism. Derivating the morphism at the identity $\delta$ yields a Lie group morphism
$$\Lf (\iota) \coloneq T_\delta \iota \colon T_\delta (\allseriesdeltamLC ) \rightarrow T_\delta (\allseriesdelta),\quad v \mapsto \iota(v) (=v).$$
Observe that the Lie bracket on $\Lf (\delta + \allseriesmLC) =T_\delta (\delta +\allseriesmLC) \cong \allseriesmLC$ coincides (pointwise) with the one of $\Lf (\delta + \allseries)$, and the latter is \eqref{compositionLiebracket}.
\endpr

Regularity of the \Frechet Lie group $\delta + \allseries$ is investigated next. For a curve $\gamma_\delta (t) = (\delta + \gamma(t)) \in \delta + \allseries$ consider the Lie type differential equation
\begin{align}\label{Lietype}
 \begin{cases}
  \dot{\gamma_\delta}(t) = \gamma_\delta(t).c(t) = c (t) + \gamma(t) \lhd c(t)  \\
  \gamma_\delta (0) = \delta,
 \end{cases}
\end{align}
where $c \colon [0,1] \rightarrow \allseries$ is a continuous curve. For every $\eta \in X^*$ observe that $(\gamma_\delta(t),\eta) = (\gamma(t),\eta)$.
Now since the coordinate functions are continuous linear, a differential equation is obtained for every word $\eta \in X^*$:
\begin{align}
 (\dot{\gamma}_\delta(t),\eta) &= (c (t),\eta) + (\gamma(t) \lhd c(t),\eta) \notag \\
 &= (c(t),\eta) + \sum_{\rho \in X^*} (\gamma(t),\rho) (\rho \lhd c(t), \eta) %\notag\\
 = (c(t),\eta) +\sum_{1 \leq |\rho| \leq |\eta|} (\gamma(t),\rho) (\rho \lhd c(t),\eta).\label{uglyODE}
\end{align}
The computations in Example \ref{ex:preLie} have been used above, and the products of elements in $\K^m$ are taken as componentwise products.
Note now that the sum in \eqref{uglyODE} only appears if $|\eta|_{x_0} \neq 0$. Hence, if a word does not contain the letter $x_0$, then the differential equation \eqref{uglyODE} reduces to
\begin{align}
(\gamma (t),\eta) = \int_0^t (c(s),\eta) \mathrm{d}s, \quad  \forall \eta \in X^\ast, |\eta|_{x_0}=0. \label{startODE}
\end{align}
Since $(c(t),\eta)$ is a continuous $\K^m$-valued curve, one can solve the above equation for all $t \in [0,1]$. Now if $\eta$ is a word with $|\eta|_{x_0} \neq 0$, observe that all elements in \eqref{uglyODE} appearing as coefficients of evaluations of $\gamma$ are continuous $\K^m$-valued curves of the form $(c(t),\eta)$ or
\begin{align}
C_{\rho,\eta} \colon [0,1] \rightarrow \K^m,\quad t \mapsto C_{\rho,\eta}(t)\coloneq (\rho \lhd c(t),\eta). \label{matrixcoeff}
\end{align}
It is now proved via induction on the length of the words
that equation \eqref{uglyODE} admits a solution on $[0,1]$ for every word.
Note first that for any word without an $x_0$ (such as the empty word, which is the only length zero element), the statement follows directly from the integral equation \eqref{startODE}.
If $|\eta| = n >1$ assume that the statement is true for all words of lower length. If $|\eta|_{x_0} =0$, the statement follows again from \eqref{startODE}. To obtain solutions for the words of length $n$ containing $x_0$, pick an enumeration $(\eta_i)_{i \in I_n}$ of words of length $n$. Using the enumeration and \eqref{matrixcoeff}, define
\begin{align*}
 \bm{v}_n (t) &\coloneq \begin{bmatrix} (\gamma(t),\eta_1) \\ (\gamma(t),\eta_2) \\ \vdots \\ (\gamma(t),\eta_{|I_n|})\end{bmatrix}, \quad C_n (t) \coloneq \begin{pmatrix}
     0 & C_{\eta_1,\eta_2}(t) &  \cdots & C_{\eta_1,\eta_{|I_n|}}(t)\\
     C_{\eta_2,\eta_1} (t) &0&  \ddots & \vdots \\
     \vdots & \ddots & \ddots & C_{\eta_{|I_n|-1}, \eta_{|I_n|}} \\
     C_{\eta_{|I_n|},\eta_1} (t) & \cdots & C_{\eta_{|I_n|}, \eta_{|I_n|-1}} & 0
    \end{pmatrix},\\ \bm{b}_n (t) &\coloneq  \sum_{|\rho| < n} \begin{bmatrix}
    (\gamma(t),\rho) (\rho \lhd c(t),\eta_1)\\
     (\gamma(t),\rho) (\rho \lhd c(t),\eta_2)\\
     \vdots \\
 (\gamma(t),\rho) (\rho \lhd c(t),\eta_{|I_n|})
   \end{bmatrix} 
\end{align*}
Then \eqref{uglyODE} together with the observation that $(\eta \lhd c , \eta)=0$ give rise to the following inhomogeneous system of linear differential equations on $(\K^m)^{|I_n|}$:
\begin{align}\label{lin:sys}
 \dot{\bm{v}}_n (t) = C_n (t) \bm{v}_n(t) + \bm{b}_n (t), \qquad t\in [0,1],
\end{align}
Now by the induction hypothesis the inhomogeneity $\bm{b}_n$ in \eqref{lin:sys} is already completely determined by the previous computations. Furthermore, the coefficient matrix $C_n$ is determined by $c$ and thus continuous in $t$. Hence, one can solve the system \eqref{lin:sys} and obtain a solution on $[0,1]$ (via the usual solution theory for linear differential equations on finite-dimensional spaces). This completes the induction, and thus, one can iteratively solve the inhomogeneous linear system \eqref{lin:sys} for every $n \in \mathbb{N}_0$ with a unique solution on $[0,1]$.
Following \cite[\S 6]{MR0463601} (cf.\ also \cite{MR1304355}), the solution to the Lie type equation \eqref{Lietype} is the solution to the infinite system of differential equations \eqref{startODE} and
$$\dot{\bm{v}}_n(t) = C_n (t)\bm{v}_n(t) + \bm{b}_n(t), \quad n\in \mathbb{N}_0.$$
The earlier discussion has shown that this system is lower diagonal, i.e., the right-hand side of the equation in degree $n$ depends only on the solutions up to degree $n$. One can now solve the differential equation on the \Frechet space by adapting the argument in \cite[p.~79-80]{MR0463601}: Lower diagonal systems can be solved iteratively
component-by-component, if each
solution exists on a time interval $[0, \varepsilon]$ for some fixed $\varepsilon > 0$. Choosing $\varepsilon =1$, observe that the Lie type equation \eqref{Lietype} admits a unique global solution which can be computed iteratively. Thus, the following result is evident.

\begin{propose}
 The Lie group $\delta + \allseriesm$ is $C^0$-regular.
\end{propose}

\begpr
It was seen in the discussion above that the \Frechet Lie group $\delta + \allseriesm$ is $C^0$-semiregular. However, due to \cite[Corollary D]{MR4000583} every $C^0$-semiregular Lie group modeled on a \Frechet space is already $C^0$-regular.
\endpr

Observe that one can leverage the regularity of the \Frechet Lie group in the investigation of the regularity for the Silva Lie group $\delta +\allseriesmLC$. The inclusion $\iota \colon \delta + \allseriesmLC \rightarrow \delta + \allseriesm$ is a Lie group morphism which relates the solutions of the evolution equation on the Silva and the \Frechet Lie group. Indeed, \cite[1.16]{HG15reg} shows that for a continuous curve $c \colon [0,1] \rightarrow \Lf(\delta +\allseriesmLC)=\allseriesmLC$ a solution to the evolution equation \eqref{Lietype} in $\delta + \allseriesmLC$ must satisfy
$$\iota \circ \Evol_{\delta + \allseriesmLC} (c) = \Evol_{\delta + \allseriesm} (\Lf(\iota) \circ c) = \Evol_{\delta + \allseriesm} (c),$$
where $c$ is interpreted canonically as a curve into $\Lf(\delta+\allseriesm)=\allseriesm$ via the natural inclusion. Hence, the Silva Lie group will be $C^0$-semiregular if and only if it can be proved that the solutions to the evolution equation on the \Frechet Lie group are bounded when the curve $c$ is bounded. Unfortunately, at present it is not obvious how to bound these solutions to the evolution equation, which leads to the following.
\smallskip

\textbf{Open problem}: Is the Silva Lie group $\delta + \allseriesmLC$ $C^0$-semiregular?

\begin{rem}
\begin{enumerate}
 \item Note that words which do not contain the letter $x_0$ do not yield the necessary bound for the solution of the evolution equation as the differential equation reduces to the integral equation \eqref{startODE} for these words.
 \item For words which contain the letter $x_0$, the linear system \eqref{lin:sys} governs the evolution equation. A natural Ansatz for the problem would thus be to apply a Gronwall type argument. Looking closer at the pre-Lie product, one easily sees that the top-level words (i.e., of length $n$ when dealing with length $n$-words) only yield an exponential bound in the Gronwall argument. Unfortunately, there seems to be no clear way to bound the norm of the inhomogeneity $\mathbf{b}_n$ in \eqref{lin:sys}.
 \item Observe that regularity of the Silva Lie group $\delta +\allseriesmLC$ follows almost directly once $C^0$-semiregularity is known: Having the semiregularity in place, it is assumed that the estimates will directly yield that for every curve $c$ taking values in
 $$\Lf(\delta + \allseriesm) \cap B_1^{\lVert \cdot \rVert_M}(0) = \{x \in \allseriesm \mid \lVert x \rVert_M \leq 1\}, M>0,$$
the evolution $\Evol (c)$ is contained in $B_K^{\lVert \cdot \rVert_N} (0)$, $N,K >0$ fixed (but depending on $M$). If this is true, $C^1$-regularity of $\delta+\allseriesmLC$ follows from the arguments presented in the proof of \cite[Theorem 4.3]{BS16}.
\end{enumerate}
\end{rem}

\appendix
\section{Infinite-dimensional calculus}\label{app:smooth}
In this appendix we recall some basic definitions concerning the infinite-dimensional calculus used throughout the article. For more information we refer to the presentations in \cite{hg2002a,neeb2006}.
\begin{definition}
 Let $r \in \mathbb{N}_0 \cup \{\infty\}$ and $E$, $F$ locally convex $\K$-vector spaces and $U \subseteq E$ open.
 A map $f \colon U \rightarrow F$ is called a $C^r_\K$-map if it is continuous and the iterated directional derivatives
  \begin{displaymath}
   d^kf (x,y_1,\ldots,y_k) \coloneq (D_{y_k} \cdots D_{y_1} f) (x)
  \end{displaymath}
 exist for all $k \in \mathbb{N}_{0}$ with $k \leq r$ and $y_1,\ldots,y_k \in E$ and $x \in U$,
 and the mappings $d^kf \colon U \times E^k \rightarrow F$ so obtained are continuous.
 If $f$ is $C^\infty_\R$, it is called \emph{smooth}.
 If $f$ is $C^\infty_\C$, it is said to be \emph{complex analytic} or \emph{holomorphic} and that $f$ is of class $C^\omega_\C$.\footnote{Recall from \cite[Proposition 1.1.16]{dahmen2011}
 that $C^\infty_\C$ functions are locally given by series of continuous homogeneous polynomials (cf.\ \cite{BS71a,BS71b}).
 This justifies the abuse of notation.}     \label{mydefinition: analyt}
\end{definition}

\begin{definition}[Complexification of a locally convex space]
 Let $E$ be a real locally convex topological vector space. Endow the locally convex product $E_\C \coloneq E \times E$ with the following operation
 \begin{displaymath}
  (x+iy).(u,v) \coloneq (xu-yv, xv+yu), \quad \forall x,y \in \R,\; u,v \in E.
 \end{displaymath}
 The complex vector space $E_\C$ is called the \emph{complexification} of $E$. Identify $E$ with the closed real subspace $E\times \{0\}$ of $E_\C$.
\end{definition}

\begin{definition}							\label{mydefinition: real_analytic}							
 Let $E$, $F$ be real locally convex spaces and $f \colon U \rightarrow F$ defined on an open subset $U$.
 $f$ is called \emph{real analytic} (or $C^\omega_\R$) if $f$ extends to a $C^\infty_\C$-map $\tilde{f}\colon \tilde{U} \rightarrow F_\C$ on an open neighborhood $\tilde{U}$ of $U$ in the complexification $E_\C$.
\end{definition}

For $r \in \mathbb{N}_0 \cup \{\infty, \omega\}$, being of class $C^r_\K$ is a local condition,
i.e.\ if $f|_{U_\alpha}$ is $C^r_\K$ for every member of an open cover $(U_\alpha)_{\alpha}$ of its domain,
then $f$  is $C^r_\K$. (See \cite[pp. 51-52]{hg2002a} for the case of $C^\omega_\R$, the other cases are clear by definition.)
In addition, the composition of $C^r_\K$-maps (if possible) is again a $C^r_\K$-map (cf. \cite[Propositions 2.7 and 2.9]{hg2002a}).

\begin{definition}[$C^r_\K$-Manifolds and $C^r_\K$-mappings between them]
 For $r \in \mathbb{N}_0 \cup \{\infty, \omega\}$, manifolds modeled on a fixed locally convex space can be defined as usual.
 Direct products of locally convex manifolds, tangent spaces and tangent bundles as well as $C^r_\K$-maps between manifolds may be defined as in the finite-dimensional setting.

For $C^r_\K$-manifolds $M,N$ the notation $C^r_\K(M,N)$ denotes the set of all $C^r_\K$-maps from $M$ to $N$.
Furthermore, for $s \in \{\infty,\omega\}$ define the \emph{locally convex $C^s_\K$-Lie groups} as groups with a $C^s_\K$-manifold structure turning the group operations into $C^s_\K$-maps.
\end{definition}

The following lemma seems to be part of the mathematical folklore, a proof can be found in \cite[Lemma A.3]{BS16}.
\begin{lemma}\label{lem:folklore}
 Let $U$ be an open subset of a complex locally convex space $E$ and $F$ be a complex locally convex space which is sequentially complete. Consider a set $\Lambda \subseteq L(F,\C)$ of complex linear functionals which separates the points on $F$.\footnote{That is, for each $x\in F$ there is a $\lambda \in \Lambda$ with $\lambda (x) \neq 0$.} If a map $f \colon U \rightarrow F$ is continuous and
  \begin{displaymath}
   \lambda \circ f \colon U \rightarrow \C
  \end{displaymath}
 is complex analytic for each $\lambda \in \Lambda$, then $f$ is complex analytic.
\end{lemma}
\bibliographystyle{new}
\bibliography{fliess_lit}

\end{document}